\documentclass[10pt]{amsart}
\usepackage{amsmath,amsthm,enumerate,graphics,pstricks,amsfonts,latexsym,amsopn,verbatim,amscd,amssymb}
\usepackage{color}
\usepackage{psfrag}
\usepackage[all]{xy}
\usepackage{tikz}
\usepackage{pgfplots}
\usepackage{caption}
\usepackage{geometry}

\theoremstyle{plain}
\newtheorem{thm}{Theorem}[section]

\newtheorem{lem}[thm]{Lemma}
\newtheorem{cor}[thm]{Corollary}

\newtheorem{problem}{Theorem}[section] 
\newtheorem{prob}[problem]{Problem}

\theoremstyle{definition}

\DeclareMathOperator{\Spec}{Spec}
\DeclareMathOperator{\L-Spec}{L-spec}
\DeclareMathOperator{\Q-Spec}{Q-spec}

\DeclareMathOperator{\E}{E}
\DeclareMathOperator{\LE}{LE}
\DeclareMathOperator{\SE}{LE^+} 
\DeclareMathOperator{\DE}{E_{D}}
\DeclareMathOperator{\DLE}{E_{DL}}
\DeclareMathOperator{\DQE}{E_{DQ}}

\begin{document}
\title[Various spectral aspects of  NCCC-graphs of certain finite non-abelian groups]{Various spectral aspects of  NCCC-graphs of certain finite non-abelian groups}

\author[R. Chakraborty, F. E. Jannat and R. K. Nath]{Rishabh Chakraborty, Firdous Ee Jannat and Rajat  Kanti Nath$^*$}
\address{Rishabh Chakraborty, Department of Mathematical Science, Tezpur  University, Napaam -784028, Sonitpur, Assam, India.}

\email{rchakraborty2101@gmail.com}

\address{Firdous Ee Jannat, Department of Mathematical Science, Tezpur  University, Napaam -784028, Sonitpur, Assam, India.}

\email{firdusej@gmail.com}
\address{Rajat  Kanti Nath, Department of Mathematical Science, Tezpur  University, Napaam -784028, Sonitpur, Assam, India.}

\email{rajatkantinath@yahoo.com}
\thanks{$^*$Corresponding author}

\begin{abstract}
  Let ${G}$ be a finite non-abelian group. The non-commuting conjugacy class graph (abbreviated as NCCC-graph) of $G$ is a simple undirected graph whose vertex set is the set of conjugacy classes of non-central elements of $G$ and two vertices $x^G$ and $y^G$ are adjacent to each other if $x'$ and $y'$ does not commute for all $x'\in x^G$ and $y'\in y^G$, where $x^G$ is the conjugacy class of $x \in G$. In this paper, we compute the spectrum, Laplacian spectrum, signless Laplacian spectrum and corresponding energies of NCCC-graphs of certain families of finite non-abelian groups. We determine whether these graphs are integral, L-integral and Q-integral.   Further,  we compare energy, Laplacian energy and signless Laplacian energy; and determine whether these graphs are borderenergetic, L-borderenergetic, Q-borderenergetic, hyperenergetic, L-hyperenergetic or Q-hyperenergetic. 
  %Additionally, we derive several  conditions (i.e. ) for these graphs. We also provide a graphical comparison of the different energies of the NCCC-graphs of  aforementioned group.   
\end{abstract}

\thanks{ }

\maketitle

%\section{Various spectra and energies of $\mathcal{E}_G$}
\section{Introduction}
Let $\mathcal{G}$ be a simple graph with vertex set $v(\mathcal{G}) = \{v_1, v_2,\ldots, v_n\}$. If $A(\mathcal{G})$ and $D(\mathcal{G})$ are the adjacency and degree matrix of $\mathcal{G}$ then the Laplacian matrix and signless Laplacian matrix of $\mathcal{G}$, denoted by $L(\mathcal{G})$ and $Q(\mathcal{G})$ respectively, are given by
\[L(\mathcal{G}) := D(\mathcal{G}) - A(\mathcal{G}) \quad \text{ and }\quad Q(\mathcal{G}) := D(\mathcal{G}) + A(\mathcal{G}).\] 
The spectrum (denoted by $\Spec(\mathcal{G})$), Laplacian spectrum (denoted by $\L-Spec(\mathcal{G})$) and signless Laplacian spectrum (denoted by $\Q-Spec(\mathcal{G})$) of $\mathcal{G}$ are sets of eigenvalues of $A(\mathcal{G})$, $L(\mathcal{G})$ and $Q(\mathcal{G})$
%adjacency matrix, Laplacian matrix and signless Laplacian matrix of $\mathcal{G}$ 
with multiplicity respectively. We write $\Spec(\mathcal{G}) = \{[\alpha_1]^{a_1}, [\alpha_2]^{a_2}, [\alpha_3]^{a_3},\ldots, [\alpha_j]^{a_j}\}$, \quad $\L-Spec(\mathcal{G}) = \{[\beta_1]^{b_1}, [\beta_2]^{b_2}, [\beta_3]^{b_3} ,\ldots, [\beta_k]^{b_k}\}$ and $\Q-Spec(\mathcal{G}) = \{[\gamma_1]^{c_1}, [\gamma_2]^{c_2}, [\gamma_3]^{c_3},\ldots, [\gamma_l]^{c_l}\}$, where $[x_i]^{y_i}$ means the algebraic multiplicity of the eigenvalue $x_i$ is $y_i$. A graph $\mathcal{G}$ is called integral, L-integral and Q-integral if $\Spec(\mathcal{G})$, $\L-Spec(\mathcal{G})$ and $\Q-Spec(\mathcal{G})$ contain only integers. The notion of integral graph was introduced by Harary and Schwenk \cite{B10} in 1974. 
%Ahmadi et al. noted that integral graphs have some interests for designing the network topology of perfect state transfer networks, see \cite{E1} and the references there in. After that 
Grone and Merris \cite{Grone-94} in 1994 and Simic and Stanic \cite{Simic-08} in 2008 introduced the notions of L-integral and Q-integral graphs respectively. %Again $\mathcal{G}$ is called super integral if $G$ is integral, L-integral and Q-integral. Then 

The energy, Laplacian energy and signless Laplacian energy of $\mathcal{G}$, denoted by $\E(\mathcal{G})$, $\LE(\mathcal{G})$ and $\SE(\mathcal{G})$ are defined as
$$\E(\mathcal{G}) := \sum\limits_{\alpha\in \Spec({\mathcal{G}})}|\alpha|, \qquad\qquad \LE(\mathcal{G}) := \sum\limits_{\beta\in \Q-Spec(\mathcal{G})}\left|\beta- \Delta(\mathcal{G})\right|$$ and $$\SE(\mathcal{G}) := \sum\limits_{\gamma\in \Q-Spec(\mathcal{G})}\left|\gamma-\Delta(\mathcal{G})\right|,$$
where $\Delta(\mathcal{G}) = \frac{|2e(\Gamma_G)|}{|v(\Gamma_G)|}$.
In 1978, Gutman \cite{Gutman78} introduced  the notion of $\E({\mathcal{G}})$; in 2006, Gutman  and Zhou \cite{zhou} introduced  the notion of $\LE({\mathcal{G}})$; and in 2008, Gutman et al. \cite{GAVBR} introduced  the notion of $\SE({\mathcal{G}})$. 
%More about these energies can be found in (Cite some papers). 
As noted in \cite{GuFu-2019}, these energies have applications in crystallography, macromolecular theory, protein sequence analysis and comparison, network analysis, satellite communication, image processing etc. 
%(see \cite{GuFu-2019} and the references therein for details). 
Note that 
\begin{equation}\label{Kn}
	\E(K_n) = \LE(K_n) = \SE(K_n) = 2(n - 1),
\end{equation}
where $K_n$ is a complete graph of order $n$. A graph $\mathcal{G}$ is called hyperenergetic if $\E(\mathcal{G}) > \E(K_{|v(\mathcal{G})|})$ and borderenergetic if $\E(\mathcal{G}) = \E(K_{|v(\mathcal{G})|})$; L-hyperenergetic if $\LE(\mathcal{G}) > \LE(K_{|v(\mathcal{G})|})$ and L-borderenergetic if $\LE(\mathcal{G}) = \LE(K_{|v(\mathcal{G})|})$ and Q-hyperenergetic if $\SE(\mathcal{G}) > \SE(K_{|v(\mathcal{G})|})$ and Q-borderenergetic if $\SE(\mathcal{G})$ $ = \SE(K_{|v(\mathcal{G})|})$.
It is worth mentioning that the concept of hyperenergetic graph was introduced by Walikar et al. \cite{Walikar} and Gutman  \cite{Gutman1999}, independently in 1999 and  L-hyperenergetic and Q-hyperenergetic graphs   were introduced in \cite{FSN-2020}. The concepts of borderenergetic, L-borderenergetic and Q-borderenergetic graphs were introduced by Gong et al. \cite{GLXGF-2015}, Tura \cite{Tura-2017} and Tao et al. \cite{TH-2018} in the years 2015, 2017 and 2018 respectively.

Let $G$ be a finite non-abelian group with center $Z(G)$. We write $x^G$ to denote the conjugacy class of $x \in G$.   The non-commuting conjugacy class graph (abbreviated as NCCC-graph) of $G$, denoted by $\Gamma_G$,  is a simple undirected graph whose vertex set is  $\{x^G: x \in G \setminus Z(G)\}$ and two vertices $x^G$ and $y^G$ are adjacent to each other if   $x'y' \ne y'x'$  for all $x'\in x^G$ and $y'\in y^G$. The complement of this graph is known as  commuting conjugacy class graph (in short CCC-graph) of $G$ and it is denoted by $CCC(G)$. The study of CCC-graphs of groups was initiated by Herzog, Longobardi, and Maj \cite{HLM} in 2009. In 2016, Mohammadian et al. \cite{MEFW-2016} characterized finite groups whose CCC-graphs are triangle-free. Subsequently,  structures of CCC-graphs of several families of finite non-abelian groups have been obtained in  \cite{SA-2020, SA-CA-2020, Salah-2020}. Later on,  Bhowal and Nath \cite{BN-2021, BN-2022, BN-2024} have computed spectrum, Laplacian spectrum, signless Laplacian spectrum and their corresponding energies. Bhowal and Nath  also obtained various groups such that CCC-graphs are hyperenergetic, L-hyperenergetic and Q-hyperenergetic.  Recently, Jannat and Nath \cite{JN-2025} have computed various common neighbourhood spectra and common neighbourhood energies of  CCC-graphs of several classes of finite non-abelian groups and obtained finite non-abelian groups such that their CCC-graphs are common neighbourhood hyperenergetic, common neighbourhood Laplacian hyperenergetic and common neighbourhood signless Laplacian hyperenergetic. All these results can  be found in the survey article \cite{CaJa-2024}. It is worth mentioning that NCCC-graphs of groups were not much studied. In our knowledge, 
\cite{JaNa7-2024} is the only paper where Jannat and Nath studied several distance spectral properties and wiener index of NCCC-graphs.

In Section 1 of this paper, we first compute spectrum, Laplacian spectrum,  signless Laplacian spectrum and their corresponding energies 
%(i.e., energy, Laplacian energy and signless Laplacian energy) 
of $\Gamma_G$ when $G$ is the group such that $\frac{G}{Z(G)}$ is isomorphic to $\mathbb{Z}_p \times \mathbb{Z}_p$ (for any prime $p$) and the dihedral group $D_{2m}$ (for any integer $m \geq 3$). As a consequence, we  compute the above mentioned graph parameters for NCCC-graphs of the dihedral group $D_{2m} = \langle x, y \mid x^m = y^2 = 1, yxy^{-1} = x^{-1} \rangle$ (for $m \geq 3$), the dicyclic group $T_{4m} = \langle x, y \mid x^{2m} = 1, x^m = y^2, y^{-1}xy = x^{-1} \rangle$ (for $m \geq 2$), the semidihedral group $SD_{8m} = \langle x, y \mid x^{4m} = y^2 = 1, yxy = x^{2m-1} \rangle$ (for $m \geq 2$) and the groups $U_{(n,m)} = \langle x, y \mid x^{2n} = y^m = 1, x^{-1}yx = y^{-1} \rangle$ (for $m \geq 3$ and $n \geq 2$),  $U_{6m} =\langle x, y : x^{2m} = y^3 = 1, x^{-1}yx = y^{-1} \rangle$ (for $m \geq 2$) and $V_{8m} = \langle x, y : x^{2m} = y^4 = 1, yx = x^{-1}y^{-1}, y^{-1}x = x^{-1}y \rangle$ (for $m \geq 2$).

 In the 
subsequent section, we compare these energies and identify conditions under which these graphs are borderenergetic, L-borderenergetic, Q-borderenergetic, hyperenergetic, L-hyperenergetic or Q-hyperenergetic. Additionally, we derive several   conditions such that these graphs  are integral, L-integral and Q-integral. We also provide a graphical comparison of the energy, Laplacian energy and signless Laplacian energy of NCCC-graphs of the aforementioned group.

\section{Various spectra and energies of NCCC-graphs}
Let $P(\mathcal{G}, x)$, $P_L(\mathcal{G}, x)$ and $P_Q(\mathcal{G}, x)$ be the characteristic polynomials of  $A(\mathcal{G})$, $L(\mathcal{G})$ and $Q(\mathcal{G})$ respectively. We write $K_{a_1 \cdot p_1 , a_2 \cdot p_2 }$ to denote the complete $r$-partite graph $K_{\underbrace{p_1,  \dots, p_1}_{a_1\text{-times}} \,, \, \underbrace{p_2,  \dots, p_2}_{a_2\text{-times}}}$, where $r = a_1 + a_2$.
The following lemma that follows from \cite[Theorem  1.6]{SN-2023} is useful in  computing spectrum, L-spectrum and Q-spectrum of the NCCC-graphs of groups.
\begin{lem}\label{Kapbq-poly}
	If  $\mathcal{G} = K_{a_1 \cdot p_1 , a_2 \cdot p_2 }$ on $n=a_1  p_1 + a_2 p_2$ vertices and $r = a_1 + a_2$ then
	\begin{enumerate}
		\item[\rm{(a)}] $P(\mathcal{G}, x) = x^{n - r}(x+p_1)^{a_1-1}(x+p_2)^{a_2-1}\,(x^2 \, + \, (p_1(1 - a_1) + p_2(1 - a_2))x + p_1p_2(1 - a_1 - a_2))$.\vspace{0.2 cm}
		\item[\rm{(b)}] $P_L(\mathcal{G}, x)= x(x - (a_1p_1 + p_2(a_2 - 1)))^{a_2(p_2-1)} (x - (a_2p_2 + p_1(a_1 - 1)))^{a_1(p_1-1)} (x - n)^{a_1 + a_2 - 1}$.\vspace{0.2 cm}
		\item[\rm{(c)}] $P_Q(\mathcal{G}, x)  =\prod_{i=1}^{2}(x - (n - p_i))^{a_i(p_i-1)} \prod_{i=1}^{2}(x - (n - 2p_i))^{a_i - 1}(x^2 - x(2n + p_1(a_1 - 2) + p_2(a_2 - 2))+ n^2 - 2n(p_1 + p_2) + n(a_1p_1 + a_2p_2) + 2p_1p_2(2 - a_1 - a_2))$. 
	\end{enumerate}	
\end{lem}
Putting $a_1 = a, p_1 = b$ and $a_2 =  p_2 = 0$ in Lemma \ref{Kapbq-poly}, we get the following result.
\begin{lem}\label{Kab-poly}
	If  $\mathcal{G} = K_{a \cdot b}$ then
	\begin{enumerate}
		\item[\rm{(a)}] $P(\mathcal{G}, x)=x^{a(b-1)}(x+b)^{a-1}\left(x- b(a - 1)\right)$.\vspace{0.2 cm} 
		\item[\rm{(b)}] $P_L(\mathcal{G}, x)= x(x - b(a-1))^{a(b-1)} (x - ab)^{a - 1}$.\vspace{0.2 cm}
		\item[\rm{(c)}] $P_Q(\mathcal{G}, x)=  (x- b(a - 1))^{a(b-1)}(x- b(a -2))^{a-1}(x-2b(a - 1))$. 
	\end{enumerate}
\end{lem}
\subsection{Groups $G$ such that $\frac{G}{Z(G)} \cong \mathbb{Z}_p \times \mathbb{Z}_p$}
In this subsection, we consider NCCC-graphs of the groups whose central quotients are isomorphic to the group  $\mathbb{Z}_p \times \mathbb{Z}_p$ for  any prime $p$.
\begin{thm}\label{Energy-NCCC1}
Let $G$ be a finite non-abelian group  such that $\frac{G}{Z(G)} \cong \mathbb{Z}_p \times \mathbb{Z}_p$, where $p$ is any prime. Then

\noindent $\Spec(\Gamma_G) = \left\{ [0]^{(p+1)(n-1)}, [-n]^{p}, [np]^1\right\}$, $\L-Spec(\Gamma_G) = \left\{[0]^1, [np]^{(p+1)(n-1)}, [(p+1)n]^{p} \right\}$, $\Q-Spec(\Gamma_G)$ $ = \left\{[np]^{(p+1)(n-1)}, [n(p-1)]^{p}, [2np]^1\right\}$
and  $\E(\Gamma_G) = \LE(\Gamma_G) = \SE(\Gamma_G) = 2np$, 
where $|Z(G)| = z$ and $n = \frac{(p-1)z}{p}$.

%\[
%\Spec(\Gamma_G) = \left\{ (0)^{(p+1)(n-1)}, (-n)^{p}, (np)^1\right\}, 
%\]
%\[
%\L-Spec(\Gamma_G) = \left\{(0)^1, (np)^{(p+1)(n-1)}, ((p+1)n)^{p} \right\},  
%\]
%\[
%\Q-Spec(\Gamma_G) = \left\{(np)^{(p+1)(n-1)}, (n(p-1))^{p}, (2np)^1\right\},
%\]
%where $n = \frac{(p-1)\left| Z\right|}{p}$. 
\end{thm}

\begin{proof}
By \cite[Theorem 3.1]{SA-2020},	we have $CCC(G) \cong (p + 1)K_n$, where $n = \frac{(p-1)z}{p}$.
Therefore, $\Gamma_G \cong K_{(p + 1) \cdot n}$. Since $|v(\Gamma_G)|=n(p+1)$,
% $a = p + 1$ and $b = n$. Therefore, 
by  Lemma \ref{Kab-poly}(a), we get 
\begin{align*}
	P(\Gamma_G, x) 
	%&= x^{(p + 1)(n - 1)}(x + n)^{p + 1 - 1}\left(x- n(p + 1 - 1)\right)\\
	  			   &= x^{(p+1)(n-1)}(x+n)^{p}\left(x- np\right).
\end{align*}
Hence, $\Spec(\Gamma_G)= \left\{ [0]^{(p+1)(n-1)}, [-n]^{p}, [np]^1\right\}$ and so
\[
\E(\Gamma_G) = (p+1)(n-1) \times 0 + p \times |-n| + |np| = 2np.
\]

 By using Lemma \ref{Kab-poly}(b), we get 
\begin{align*}
	P_L(\Gamma_G, x) 
	%&= x(x - n(p + 1 - 1))^{(p + 1)(n - 1)} (x - (p + 1)n)^{p + 1 - 1}\\
					 &= x(x - np)^{(p+1)(n-1)} (x - (p + 1)n)^{p}.
\end{align*}
Therefore, $\L-Spec(\Gamma_G) = \left\{[0]^1, [np]^{(p+1)(n-1)}, [(p+1)n]^{p} \right\}$. 

We have $|2e(\Gamma_G)| = p(p+1)n^2$ and so $\Delta(\Gamma_G) = \frac{|2e(\Gamma_G)|}{|v(\Gamma_G)|} = pn$. Therefore,
\small{
\begin{align*}
	\LE(\Gamma_G) &= 1 \times |0 - \Delta(\Gamma_G)| + (p+1)(n-1) \times |np - \Delta(\Gamma_G)| + p \times |(p+1)n - \Delta(\Gamma_G)|\\
	&= 1 \times |-pn| + (p+1)(n-1) \times |0| +  p \times |n| = 2np.
\end{align*}
}
%Hence, $\LE(\Gamma_G) = 2np$.\\
%\begin{align*} 
% 	 {|2e(\Gamma_G)|} &= 1 \times 0 + (p+1)(n-1) \times np + p \times (p+1)n.\\
% 	                         &= p(p+1)n^2.
%\end{align*} 
%and so $\frac{|2e(\Gamma_G)|}{|v(\Gamma_G)|} = pn$.\\
%Hence,
%\[
%\left|0 - \frac{2|e(\Gamma_G)|}{|v(\Gamma_G)|} \right| =  \left| pn \right| = pn,
%\]
%\begin{align*}
%	\left|\frac{2|e(\Gamma_G)|}{|v(\Gamma_G)|} - pn\right| &=  \left| pn -  pn\right| = 0 ,
%\end{align*}
%and
%\begin{align*}
%	\left|\frac{2|e(\Gamma_G)|}{|v(\Gamma_G)|}  - (p+1)n\right| &=  \left| pn -  (p+1)n\right| = n.
%\end{align*}
%Therefore,
%\begin{align*}
%	LE(\Gamma_G) &= 1 \times pn + (p+1)(n-1) \times 0 + p \times n.\\
%	&= 2np.
%\end{align*}
 By Lemma \ref{Kab-poly}(c), we get 
\begin{align*}
	P_Q(\Gamma_G, x) 
	%&= (x- n(p + 1 - 1))^{(p + 1)(n - 1)}(x - n(p + 1  - 2))^{p + 1 - 1}(x - 2n(p + 1 - 1))\\
	 				 &= (x- np)^{(p+1)(n-1)}(x- n(p - 1))^{p}(x-2np).
\end{align*}
Therefore, $\Q-Spec(\Gamma_G) = \left\{[np]^{(p+1)(n-1)}, [n(p-1)]^{p}, [2np]^1\right\}$ and so
\small{
\begin{align*}
	\SE(\Gamma_G) &= (p + 1)(n - 1) \times |np - \Delta(\Gamma_G)| + p \times | n(p - 1) - \Delta(\Gamma_G)| + 1 \times | 2np - \Delta(\Gamma_G)|\\
	& =  (p + 1)(n - 1) \times |0| + p \times | -n | + 1 \times | pn | = 2np.
\end{align*}
}
%Therefore $\SE(\Gamma_G) = 2np$. Hence the result follows.
This completes the proof.
\end{proof}

\begin{cor}
	Let $G$ be a non-abelian group of order $p^m$ and $|Z(G)| = p^{m-2}$. Then 
	%the Spec, L-spec, Q-spec and their respective energies of NCCC-graph of $G$ are given by 

 $\Spec(\Gamma_G) = \left\{ [0]^{(p+1)\left((p-1)p^{m-3}-1\right)}, [-(p-1)p^{m-3}]^{p}, [(p-1)p^{m-2}]^1\right\}$,
	
 $\L-Spec(\Gamma_G) = \left\{[0]^1, [(p-1)p^{m-2}]^{(p+1)\left((p-1)p^{m-3}-1\right)}, [(p^2 - 1)p^{m-3}]^{p} \right\}$,

 \small{$\Q-Spec(\Gamma_G) = \left\{[(p-1)p^{m-2}]^{(p+1)\left((p-1)p^{m-3}-1\right)}, [(p-1)^2p^{m-3}]^{p}, [2(p-1)p^{m-2}]^1\right\}$.}

\noindent and  $\E(\Gamma_G) = \LE(\Gamma_G) = \SE(\Gamma_G) = 2(p-1)p^{m-2}$.
\end{cor}

\begin{proof}
Here $\left|\frac{G}{Z(G)} \right| = p^2$. Since $G$ is non-abelian, $\frac{G}{Z(G)}$ is non-cyclic. Therefore,
 $\frac{G}{Z(G)} \cong \mathbb{Z}_p \times \mathbb{Z}_p$. Hence, the result follows from Theorem \ref{Energy-NCCC1}. 
\end{proof}

\subsection{Groups $G$ such that $\frac{G}{Z(G)} \cong D_{2m}$}
In this subsection, we consider NCCC-graphs of the groups whose central quotients are isomorphic to the group  $D_{2m}$.
\begin{thm}\label{Energy-NCCC2}
	Let $G$ be a finite group and $\frac{G}{Z(G)} \cong D_{2m}$, where $m \geq 2$. Suppose that  $|Z(G)| = z$. Then we have the following.
	% and $\left| Z\right|$ = $z$. If $\frac{G}{Z}$ is isomorphic to   $D_{2m}$, then the Spec, L-spec, Q-spec and their respective energies of NCCC-graph of $G$ are given by 
	\begin{enumerate}
		\item For $m$ is even\\
		\rm{(i)} $\Spec(\Gamma_G) = \left\{ [0]^{\frac{z(m + 1)}{2} - 3} , \left[ \frac{-z}{2} \right]^1 ,\left[ \frac{z\left(1 + \sqrt{8m - 7)}\right)}{4}\right]^1, \left[ \frac{z\left(1 - \sqrt{8m - 7}\right)}{4}\right]^1 \right\}$
%1, 2, 4, 7, 11, 16, 22, 29, 37, 46, 56, 67, 79, 92, 106, 121, 137, 154, 172, 191 upto 200		 
		and 
		\[\E(\Gamma_G) = \frac{z\left(1 + \sqrt{8m - 7}\right)}{2}.\]
		\rm{(ii)} $\L-Spec(\Gamma_G) = \left\{ [0]^1, \left [z\right ]^{\frac{(m-1)z}{2} - 1}, \left [\frac{zm}{2} \right ]^{z - 2}, \left [ \frac{z(m+1)}{2}\right ]^2 \right\}$ and  $\LE(\Gamma_G) =  \frac{z^2(m-1)(m-2)}{(m+1)} + 2z .$ 
		\rm{(iii)} $\Q-Spec(\Gamma_G) = \left\{\left [z\right ]^{\frac{(m-1)z}{2} - 1},  \left [\frac{zm}{2} \right ]^{z - 2}, \left[\frac{z(m - 1)}{2} \right]^1,\right. \\ \left. \qquad\qquad\qquad\qquad\qquad\qquad\qquad\left[\frac{z\left(m + 3 + \sqrt{(m - 1)(m + 7)}\right)}{4} \right]^1, \left[\frac{z\left(m + 3 - \sqrt{(m - 1)(m + 7)}\right)}{4} \right]^1 \right\}$ and 
		\[\SE(\Gamma_G) = \frac{z^2(m-1)(m-2)}{m+1} +\frac{z(m-1)\left(\sqrt{1+\frac{8}{m-1}}-1\right)}{2}.\]
		%\[\SE(\Gamma_G) = \frac{z^2(m-1)(m-2)}{m+1} +\frac{z\left(\sqrt{(m-1)(m+7)}-(m-1)\right)}{2}.\]
	\item For $m$ is odd\\
	\rm{(i)} $\Spec(\Gamma_G) = \left\{ [0]^{\frac{z(m + 1)}{2} - 2} ,\left [z\sqrt{\frac{m - 1}{2}}\right ]^1, \left [-z\sqrt{\frac{m - 1}{2}}\right ]^1\right\}$ and
	$\E(\Gamma_G) = 2z\sqrt{\frac{m - 1}{2}}.$
	
\noindent	\rm{(ii)} $\L-Spec(\Gamma_G) = \left\{ [0]^1, \left [z\right ]^{\frac{z(m-1)}{2} - 1}, \left [\frac{z(m - 1)}{2} \right ]^{z - 1}, \left [ \frac{z(m+1)}{2}\right ]^1 \right\}$ and \[\LE(\Gamma_G) = \left ( \frac{z^2(m-1)(m-3)}{(m+1)}\right ) + 2z.\]
	\rm{(iii)} $\Q-Spec(\Gamma_G) = \left\{ [0]^1, \left [z\right ]^{\frac{z(m-1)}{2} - 1}, \left [\frac{z(m - 1)}{2} \right ]^{z - 1}, \left [ \frac{z(m+1)}{2}\right ]^1 \right\}$ and \[\SE(\Gamma_G) = \left ( \frac{z^2(m-1)(m-3)}{(m+1)}\right ) + 2z.\]	
	\end{enumerate}

\end{thm}

\begin{proof}
\rm{(a)} If $m$ is even then by \cite[Theorem 1.2]{Salah-2020},	we have $CCC(G) \cong 2K_{\frac{z}{2}} \cup K_{\frac{(m - 1)z}{2}}$. Therefore, $\Gamma_G \cong K_{2 \cdot \frac{z}{2} , 1 \cdot \frac{(m - 1)z}{2}}$. 	
%	If $\frac{G}{Z(G)} \cong D_{2m}$ and $m$ is even then by \cite[Theorem 1.2]{Salah-2020},	we have 
%	$$\Gamma_G \cong K_{2 \cdot \frac{z}{2} , \frac{(m - 1)z}{2}}.$$
Here $|v(\Gamma_G)| = \frac{z(m +1)}{2}$.
% and $r =3$.\\

\rm{(i)}  Using Lemma \ref{Kapbq-poly}(a), we get
\[
P(\Gamma_G, x) = x^{\frac{z(m + 1)}{2} - 3}\left(x+\frac{z}{2}\right) \left(x^2- x\left(\frac{z}{2}\right) - \frac{z^2(m - 1)}{2}\right).
\]
%\begin{align*}
%		P(\Gamma_G, x) &=x^{\frac{z(m +1)}{2} - 3}\left(x+\frac{z}{2}\right)^{2-1}\left(x+\frac{(m - 1)z}{2}\right)^{1-1}\\
%		&\left(x^2+ \left(\frac{z}{2}(1 - 2) + \frac{(m - 1)z}{2}(1 - 1)\right)x + \frac{z}{2}\frac{(m - 1)z}{2}(1 - 2 - 1)\right)\\
%		&= x^{\frac{z(m + 1)}{2} - 3}\left(x+\frac{z}{2}\right) \left(x^2- x\left(\frac{z}{2}\right) - \frac{z^2(m - 1)}{2}\right).
%\end{align*}
Thus, $\Spec(\Gamma_G) = \left\{ [0]^{\frac{z(m + 1)}{2} - 3} , \left[ \frac{-z}{2} \right]^1 ,\left[ \frac{z\left(1 + \sqrt{1+8(m - 1)}\right)}{4}\right]^1, \left[ \frac{z\left(1 - \sqrt{1+8(m - 1)}\right)}{4}\right]^1 \right\}$ and so
\small{
\begin{align*}
	\E(\Gamma_G) &\!=\! \left({\frac{z(m + 1)}{2} - 3}\right)\! \times\! 0 \!+\!\left|-\frac{z}{2}\right| \!+\! \left|\frac{z\left(1 + \sqrt{1+8(m - 1)}\right)}{4} \right| \!+\! \left|\frac{z\left(1 - \sqrt{1+8(m - 1)}\right)}{4} \right|\\ 
	 &\!=\! \frac{z}{2} + \frac{z\left(1 + \sqrt{1+8(m - 1)}\right)}{4} + \frac{z\left(\sqrt{1+8(m - 1)}-1\right)}{4}\\
	 &\!=\! \frac{z\left(1 + \sqrt{1+8(m - 1)}\right)}{2}.
    \end{align*}
}
\rm{(ii)} Using Lemma \ref{Kapbq-poly}(b), we get
\begin{align*}
		P_L(\Gamma_G, x) 
		%&= x\left(x - \left(2\frac{z}{2} + \frac{(m - 1)z}{2}(1 - 1)\right)\right)^{1\left(\frac{(m - 1)z}{2}-1\right)}\\
		%&\left(x - \left(\frac{(m - 1)z}{2} + \frac{z}{2}(2 - 1)\right)\right)^{2(\frac{z}{2}-1)} \left(x - \frac{z(m +1)}{2}\right)^{2 + 1 - 1}\\
		&= x(x - z)^{\frac{(m-1)z}{2} - 1} \left(x - \frac{zm}{2}\right)^{z - 2} \left(x - \frac{z(m+1)}{2}\right)^2.
	\end{align*}
Therefore, $\L-Spec(\Gamma_G) = \left\{ [0]^1, \left [z\right ]^{\frac{(m-1)z}{2} - 1}, \left [\frac{zm}{2} \right ]^{z - 2}, \left [ \frac{z(m+1)}{2}\right ]^2 \right\}$. 

Now $|2e(\Gamma_G)|= \frac{z^2(2m-1)}{2}$ and so $\Delta(\Gamma_G) = \frac{|2e(\Gamma_G)|}{|v(\Gamma_G)|} =  \frac{z(2m-1)}{(m+1)}$. We have
%We have \begin{align*} 
%		|2e(\Gamma_G)|= \left ( \frac{z^2(2m-1)}{2}\right ).
%\end{align*} 
%	and so $\frac{|2e(\Gamma_G)|}{|v(\Gamma_G)|} = \left ( \frac{z(2m-1)}{(m+1)}\right )$.\\
%	Hence,
\[
\left|0 - \Delta(\Gamma_G) \right| =  \left|  \frac{-z(2m-1)}{(m+1)} \right| =  \frac{z(2m-1)}{(m+1)}, \quad
	\left|z - \Delta(\Gamma_G)\right| =  \left| \frac{-z(m - 2)}{(m+1)}\right| =  \frac{z(m - 2)}{(m+1)},
\]
\begin{align*}
	\left|\frac{zm}{2} - \Delta(\Gamma_G) \right| &=  \left|  \frac{z(m - 2)(m-1)}{2(m+1)}\right| = \frac{z(m - 2)(m-1)}{2(m+1)}
\end{align*}
and
\begin{align*}
	\left|\frac{z(m+1)}{2} - \Delta(\Gamma_G)  \right| &=  \left|  \frac{z\left((m-1)^2+2\right)}{2(m+1)}\right| =  \frac{z((m-1)^2+2)}{2(m+1)}.
\end{align*}
Therefore,
\begin{align*}
	\LE(\Gamma_G) &= 1 \times \left| \frac{-z(m - 2)}{(m+1)}\right| + \left(\frac{z(m-1)}{2} - 1\right) \times \left| \frac{-z(m - 2)}{(m+1)}\right|\\
	&+ (z-2) \times \left|  \frac{z(m - 2)(m-1)}{2(m+1)}\right| + 2 \times  \left|  \frac{z\left((m-1)^2+2\right)}{2(m+1)}\right|\\
	&= \frac{z^2(m-1)(m-2)}{(m+1)} + 2z.
\end{align*}
\rm{(iii)} Using Lemma \ref{Kapbq-poly}(c), we get	
\begin{align*}
    P_Q(\Gamma_G, x) 
    %&= \prod_{i=1}^{2}\left(x - \left( \frac{z(m+1)}{2} - p_i\right)\right)^{a_i(p_i-1)} \prod_{i=1}^{2}\left(x - \left( \frac{z(m+1)}{2} - 2p_i\right)\right)^{a_i - 1}\\
   % &\bigg(x^2- x\left(2\frac{z(m+1)}{2}- p_2\right)+ \left(\frac{z(m+1)}{2}\right)^2 - 2\frac{z(m+1)}{2}(p_1 + p_2)\\
    %&\qquad\qquad\qquad\qquad\qquad\qquad\qquad\qquad\qquad + \frac{z(m+1)}{2}(2p_1 + p_2) - 2p_1p_2\bigg)\\
	&=\left(x - z\right)^{\frac{(m-1)z}{2} - 1}\left(x-\frac{zm}{2}\right)^{z - 2} \left(x - \frac{z(m -1)}{2}\right)\biggl(x^2- x\\ 
	&\left(2\frac{z(m+1)}{2}- \frac{(m -1)z}{2}\right) +\left(\frac{z(m+1)}{2}\right)^2-2\frac{z(m+1)}{2}\left(\frac{z}{2} + \frac{(m - 1)z}{2}\right) \\
	&\qquad\qquad\qquad\qquad\qquad\qquad + \frac{z(m+1)}{2}\left(2\frac{z}{2}+ \frac{(m - 1)z}{2}\right) - 2\frac{z}{2}\frac{(m - 1)z}{2} \biggr).
   \end{align*}
	Therefore,
	\small{
		\begin{align*}
		\Q-Spec(\Gamma_G) &= \left\lbrace\left [z\right ]^{\frac{(m-1)z}{2} - 1},  \left [\frac{zm}{2} \right ]^{z - 2}, \left[\frac{z(m - 1)}{2} \right]^1,\right.\\ &\left.\left[\frac{z\left(m + 3 + \sqrt{(m-1)(m+7)}\right)}{4} \right]^1, \left[\frac{z\left(m + 3 - \sqrt{(m-1)(m+7)}\right)}{4} \right]^1 \right\rbrace.
	\end{align*}
%		We have \begin{align*} 
%		{|2e(\Gamma_G)|}= \left ( \frac{z^2(2m-1)}{2}\right ).
%	\end{align*} 
%	and so $\frac{|2e(\Gamma_G)|}{|v(\Gamma_G)|} = \left ( \frac{z(2m-1)}{(m+1)}\right )$.\\
%	Hence,
We have	
\begin{align*}
	\left|\frac{z(m-1)}{2} - \Delta(\Gamma_G)  \right| &=  \left| \frac{z((m-2)^2-3)}{2(m+1)}\right| =  \frac{z((m-2)^2-3)}{2(m+1)},
\end{align*}
\small{
\begin{align*}
	\left|\frac{z\left(m + 3 + (m - 1)\sqrt{1 + \frac{8}{m - 1}}\right)}{4}\! -\! \Delta(\Gamma_G)\right| &=  \left| \frac{z\left((m - 2)^2 + 3 + (m + 1)\sqrt{(m - 1)(m+7)}\right)}{4(m+1)} \right|\\ 
	&= \frac{z\left((m - 2)^2 + 3 + (m + 1)\sqrt{(m - 1)(m+7)}\right)}{4(m+1)}
\end{align*}
}
and
\small{
\begin{align*}
	\left|\frac{z\left(m + 3 - (m - 1)\sqrt{1 + \frac{8}{m - 1}}\right)}{4}\! -\! \Delta(\Gamma_G) \right| &=  \left| \frac{z\left((m - 2)^2 + 3 - (m + 1)\sqrt{(m - 1)(m+7)}\right)}{4(m+1)}\right|\\ 
	 &= \frac{-z\left((m - 2)^2 + 3 - (m + 1)\sqrt{(m - 1)(m+7)}\right)}{4(m+1)},
\end{align*}
}
since for $m \geq 4$, $((m - 2)^2 + 3)^2 - (m + 1)^2(m - 1)(m+7) = -8m^2(2m - 3) - 48m + 56 < 0$ and so $(m - 2)^2 + 3 < (m + 1)\sqrt{(m - 1)(m+7)}$. Therefore, 
\begin{align*}
	\SE(\Gamma_G) &= \left(\frac{(m-1)z}{2} - 1 \right) \times \left| \frac{-z(m - 2)}{(m+1)}\right| + (z-2) \times \left|  \frac{z(m - 2)(m-1)}{2(m+1)}\right|\\
	&+ 1 \times \left| \frac{z((m-2)^2-3)}{2(m+1)}\right|
	+\left| \frac{z\left((m - 2)^2 + 3 + (m + 1)\sqrt{(m - 1)(m+7)}\right)}{4(m+1)} \right|\\
	&+ \left| \frac{z\left((m - 2)^2 + 3 - (m + 1)\sqrt{(m - 1)(m+7)}\right)}{4(m+1)}\right|\\
	&= \left(\frac{z^2(m-1)(m-2)}{m+1}\right) + \left(\frac{z(m-1)\left(\sqrt{1+\frac{8}{m-1}}-1\right)}{2}\right).
	%&= \frac{z^2(m-1)(m-2)}{m+1} +\frac{z\left(\sqrt{(m-1)(m+7)}-(m-1)\right)}{2}.
\end{align*}
Hence the result follows.}
\\
\rm{(b)} If $m$ is odd then by \cite[Theorem 1.2]{Salah-2020}, we have $CCC(G) \cong K_z \cup K_{\frac{(m - 1)z}{2}}$. Therefore, $\Gamma_G \cong K_{1\cdot z , 1 \cdot\frac{(m - 1)z}{2}}$.
Here $|v(\Gamma_G)| = \frac{z(m +1)}{2}$ and $r =2$.\\
\rm{(i)} Using Lemma \ref{Kapbq-poly}(a), we get	
	\begin{align*}
		P(\Gamma_G, x) 
		%&=x^{\frac{z(m + 1)}{2} - 2}(x+z)^{1-1}(x+\frac{(m - 1)z}{2})^{1-1}\\
		%&\left(x^2+ \left( z(1 - 1) + \frac{(m - 1)z}{2}(1 - 1) \right)x + z\frac{(m - 1)z}{2}(1 - 1 - 1)\right)\\
		&= x^{\frac{z(m + 1)}{2} - 2}\left(x^2 - \frac{z^2(m - 1)}{2}\right).
	\end{align*}
	Thus, $\Spec(\Gamma_G) = \left\{ [0]^{\frac{z(m + 1)}{2} - 2} ,\left [z\sqrt{\frac{m - 1}{2}}\right ]^1, \left [-z\sqrt{\frac{m - 1}{2}}\right ]^1\right\}$ and so
	\begin{align*}
		\E(\Gamma_G) &= \left({\frac{z(m + 1)}{2} - 2}\right) \times 0 + 1 \times \left|z\sqrt{\frac{m - 1}{2}}\right| + 1 \times \left|-z\sqrt{\frac{m - 1}{2}} \right|
		%&= z\sqrt{\frac{m - 1}{2}} + z\sqrt{\frac{m - 1}{2}}\\
		= 2z\sqrt{\frac{m - 1}{2}}.
	\end{align*}
\rm{(ii)} Using Lemma \ref{Kapbq-poly}(b), we get	 
	\begin{align*}
		P_L(\Gamma_G, x) 
		%&= x\left(x - \left(z + \frac{(m - 1)z}{2}(1 - 1)\right)\right)^{(\frac{(m - 1)z}{2}-1)}\\
		%&\left(x - \left(\frac{(m - 1)z}{2} + z(1 - 1)\right)\right)^{(z-1)} \left(x - \frac{z(m + 1)}{2}\right)^{1 + 1 - 1}\\
		&= x(x - z)^{\frac{z(m-1)}{2} - 1} \left(x - \frac{z(m - 1)}{2}\right)^{z - 1} \left(x - \frac{z(m+1)}{2}\right)^1 .
\end{align*}
So, $\L-Spec(\Gamma_G) = \left\{ [0]^1, \left [z\right ]^{\frac{z(m-1)}{2} - 1}, \left [\frac{z(m - 1)}{2} \right ]^{z - 1}, \left [ \frac{z(m+1)}{2}\right ]^1 \right\}$ .\\
Now $|2e(\Gamma_G)|= \left ( z^2(m-1)\right )$ and so $\Delta(\Gamma_G) = \frac{|2e(\Gamma_G)|}{|v(\Gamma_G)|} = \frac{2z(m-1)}{(m+1)}$. We have
%We have \begin{align*} 
%We have \begin{align*} 
%	{|2e(\Gamma_G)|}= \left ( z^2(m-1)\right ).
%\end{align*} 
%and so $\frac{|2e(\Gamma_G)|}{|v(\Gamma_G)|} = \left ( \frac{2z(m-1)}{(m+1)}\right )$.\\
%Hence,
\[
\left|0 - \Delta(\Gamma_G) \right| =  \left|  \frac{-2z(m-1)}{(m+1)}\right| =   \frac{2z(m-1)}{(m+1)},
\quad
	\left|z - \Delta(\Gamma_G) \right| =  \left| \frac{-z(m - 3)}{(m+1)} \right| =  \frac{z(m - 3)}{(m+1)},
\]
\begin{align*}
	\left|\frac{z(m-1)}{2} - \Delta(\Gamma_G) \right| &=  \left| \frac{z(m - 1)(m-3)}{2(m+1)}\right| =  \frac{z(m - 1)(m-3)}{2(m+1)}
\end{align*}
and
\begin{align*}
	\left|\frac{z(m+1)}{2} - \Delta(\Gamma_G) \right| &=  \left| \frac{z\left((m-1)^2+4\right)}{2(m+1)}\right| =  \frac{z\left((m-1)^2+4\right)}{2(m+1)}.
\end{align*}

Therefore,
\begin{align*}
	\LE(\Gamma_G) &= 1 \times  \left|  \frac{-2z(m-1)}{(m+1)}\right| + \left(\frac{z(m-1)}{2} - 1\right) \times  \left| \frac{-z(m - 3)}{(m+1)} \right|\\
	&+ (z-1) \times  \left| \frac{z(m - 1)(m-3)}{2(m+1)}\right|  + 1 \times \left| \frac{z\left((m-1)^2+4\right)}{2(m+1)}\right|\\
	&= \frac{z^2(m-1)(m-3)}{(m+1)} + 2z .
\end{align*}
\rm{(iii)} Using Lemma \ref{Kapbq-poly}(c), we get	
	\begin{align*}
	P_Q(\Gamma_G, x)& 
	%= \prod_{i=1}^{2}\left(x - \left( \frac{z(m+1)}{2} - p_i\right)\right)^{a_i(p_i-1)} \prod_{i=1}^{2}\left(x - \left( \frac{z(m+1)}{2} - 2p_i\right)\right)^{a_i - 1}\\
	%&\qquad\qquad\qquad\left(x^2 - x\left(2\frac{z(m+1)}{2}- p_1 - p_2\right)+ \left(\frac{z(m+1)}{2}\right)^2 \right. \\ &\qquad\qquad\qquad\qquad\qquad\left. - 2\frac{z(m+1)}{2}(p_1 + p_2) + \frac{z(m+1)}{2}(p_1 + p_2) \right)\\
	= \left(x - z\right)^{\frac{z(m-1)}{2} - 1}\left(x - \left(\frac{z(m - 1)}{2} \right)\right)^{z - 1}\biggl(x^2 - x\left(2\frac{z(m+1)}{2}- z - \frac{z(m-1)}{2}\right)\\
	&+ \left(\frac{z(m+1)}{2}\right)^2 - 2\frac{z(m+1)}{2}\left(z + \frac{z(m-1)}{2}\right) + \frac{z(m+1)}{2}\left(z + \frac{z(m-1)}{2}\right)\biggr).
\end{align*}
Therefore,  $\Q-Spec(\Gamma_G) = \left\{ [0]^1, \left [z\right ]^{\frac{z(m-1)}{2} - 1}, \left [\frac{z(m - 1)}{2} \right ]^{z - 1}, \left [ \frac{z(m+1)}{2}\right ]^1 \right\}$.
Since $\L-Spec(\Gamma_G) = \Q-Spec(\Gamma_G)$, so  $\LE(\Gamma_G) = \SE(\Gamma_G)$.
\end{proof}

\begin{cor}\label{C1}
Consider the dihedral group $D_{2m}$, where $m \geq 3$.  
\begin{enumerate}
\item For $m$ is odd\\
		\rm{(i)} $\Spec(\Gamma_{D_{2m}}) = \left\{ [0]^{\frac{(m -3)}{2}} ,\left [\sqrt{\frac{m - 1}{2}}\right ]^1, \left [-\sqrt{\frac{m - 1}{2}}\right ]^1\right\}$ and
		$\E(\Gamma_{D_{2m}}) = 2\sqrt{\frac{m - 1}{2}}$.
		
\noindent		\rm{(ii)} $\L-Spec(\Gamma_{D_{2m}}) = \left\{ [0]^1, \left [1\right ]^{\frac{m-3}{2}}, \left [ \frac{m+1}{2}\right ]^1 \right\}$ and $\LE(\Gamma_{D_{2m}}) =  \frac{(m-1)(m-3)}{(m+1)} + 2$.
		
\noindent		\rm{(iii)} $\Q-Spec(\Gamma_{D_{2m}}) = \left\{ [0]^1, \left [1\right ]^{\frac{(m-3)}{2}}, \left [ \frac{m+1}{2}\right ]^1 \right\}$ and $\SE(\Gamma_{D_{2m}}) = \frac{(m-1)(m-3)}{(m+1)} + 2$.	
		\item For $m$ is even and $\frac{m}{2}$ is odd
		
		\noindent\rm{(i)} $\Spec(\Gamma_{D_{2m}}) = \left\{ [0]^{\frac{(m - 2)}{2}} , \left [\sqrt{m-2}\right ]^1, \left [-\sqrt{m-2}\right ]^1  \right\}$ and 
		$\E(\Gamma_{D_{2m}}) = 2\sqrt{m - 2}$.
		
		\noindent\rm{(ii)} $\L-Spec(\Gamma_{D_{2m}}) = \left\{ [0]^1, \left [2\right ]^{\frac{m-4}{2}},\left[ \frac{m-2}{2}\right]^1, \left [ \frac{m+2}{2}\right ]^1 \right\}$ and  $\LE(\Gamma_{D_{2m}}) =  \frac{2(m-2)(m-6)}{(m+2)} + 4$. 
		
		\noindent \rm{(iii)} $\Q-Spec(\Gamma_{D_{2m}}) = \left\{ [0]^1, \left [2\right ]^{\frac{(m-4)}{2}}, \left [ \frac{m-2}{2}\right ]^1, \left [ \frac{m+2}{2}\right ]^1 \right\}$ and 
		$\SE(\Gamma_{D_{2m}}) = \frac{2(m-2)(m-6)}{(m+2)} + 4$.
		\item For $m$ and $\frac{m}{2}$ are even\\
		\rm{(i)} $\Spec(\Gamma_{D_{2m}}) = \left\{ [0]^{\frac{(m - 4)}{2}} , \left[ -1 \right]^1 ,\left[ \frac{1 + \sqrt{4m - 7}}{2}\right]^1, \left[ \frac{1 - \sqrt{4m - 7}}{2}\right]^1 \right\}$
		%2, 4, 8, 14, 22, 32, 44, 58, 74, 92, 112, 134, 158, 184 upto 200
		 and  
		$\E(\Gamma_{D_{2m}}) = 1 + \sqrt{4m - 7}$.
		
		\noindent \rm{(ii)} $\L-Spec(\Gamma_{D_{2m}}) = \left\{ [0]^1, \left [2\right ]^{\frac{m-4}{2}}, \left [ \frac{m+2}{2}\right ]^2 \right\}$ and  $\LE(\Gamma_{D_{2m}}) =  \frac{2(m-2)(m-4)}{(m+2)} + 4$. 
		
		\noindent \rm{(iii)} $\Q-Spec(\Gamma_{D_{2m}}) = \left\{\left [2\right ]^{\frac{m-4}{2}}, \left[\frac{m - 2}{2} \right]^1, \left[\frac{m + 6 + \sqrt{(m - 2)(m + 14)}}{4} \right]^1,\right.   \left. \left[\frac{m + 6 - \sqrt{(m - 2)(m + 14)}}{4} \right]^1 \right\}$ and 
		\[\SE(\Gamma_{D_{2m}}) = \begin{cases}
			4, &\text{for } m = 4\\
			\frac{2(m-2)(m-4)}{(m+2)} +\frac{(m-2)\left(\sqrt{1+\frac{16}{m-2}}-1\right)}{2},&\text{for } m\neq 4.
			%\frac{2(m-2)(m-4)}{(m+2)} +\frac{\left(\sqrt{(m-2)(m+14)}-(m-2)\right)}{2},&\text{for } m\neq 4 .
			\end{cases}\]
	\end{enumerate} 
\end{cor}
\begin{proof}
	We know that $\frac{D_{2m}}{Z(D_{2m})} \cong D_{2 \times \frac{m}{2}}$ or $D_{2m}$ according as $m$ is even or odd. Therefore, by Theorem \ref{Energy-NCCC1} and  Theorem \ref{Energy-NCCC2}, we get the required result.
\end{proof}

\begin{cor}\label{C2}
	Consider the dicyclic group $T_{4m}$, where $m \geq 2$.	
	\begin{enumerate}
		\item For $m$ is even\\
		\rm{(i)} $\Spec(\Gamma_{T_{4m}}) = \left\{ [0]^{m - 2} , \left[ -1 \right]^1 ,\left[ \frac{1 + \sqrt{8m - 7}}{2}\right]^1, \left[ \frac{1 - \sqrt{8m - 7}}{2}\right]^1 \right\}$ and 
		$\E(\Gamma_{T_{4m}}) = 1 + \sqrt{8m - 7}.$
		\rm{(ii)} $\L-Spec(\Gamma_{T_{4m}}) = \left\{ [0]^1, \left [2\right ]^{m- 2}, \left [ m + 1\right ]^2 \right\}$ and  $\LE(\Gamma_{T_{4m}}) =  \frac{4(m-1)(m-2)}{(m+1)} + 4.$ 
		 
		\noindent\rm{(iii)} $\Q-Spec(\Gamma_{T_{4m}}) = \left\{\left [2\right ]^{m - 2}, \left[m - 1 \right]^1, \left[\frac{m + 3 + \sqrt{(m-1)(m+7)}}{2} \right]^1,\right.   \left.\left[\frac{m + 3 - \sqrt{(m-1)(m+7)}}{2} \right]^1 \right\}$ and 
		\[\SE(\Gamma_{T_{4m}}) = \begin{cases}
		    4, &\text{for } m = 2\\
		 \frac{4(m-1)(m-2)}{(m+1)} +(m-1)\left(\sqrt{1+\frac{8}{m-1}}-1\right),&\text{for } m\neq 2.
		  %\frac{4(m-1)(m-2)}{(m+1)} +\left(\sqrt{(m-1)(m+7)}-(m-1)\right),&\text{for } m\neq 2.  
		\end{cases}\]
		\item For $m$ is odd\\
		\rm{(i)} $\Spec(\Gamma_{T_{4m}}) = \left\{ [0]^{m - 1} ,\left [2\sqrt{\frac{m - 1}{2}}\right ]^1, \left [-2\sqrt{\frac{m - 1}{2}}\right ]^1\right\}$ and
		$\E(\Gamma_{T_{4m}}) = 4\sqrt{\frac{m - 1}{2}}.$
		
		\noindent\rm{(ii)} $\L-Spec(\Gamma_{T_{4m}}) = \left\{ [0]^1, \left [2\right ]^{m - 2}, \left [ m-1\right ]^1, \left [ m+1\right ]^1 \right\}$ and $\LE(\Gamma_{T_{4m}}) =  \frac{4(m-1)(m-3)}{(m+1)} + 4.$
		
		\noindent \rm{(iii)} $\Q-Spec(\Gamma_{T_{4m}}) = \left\{ [0]^1, \left [2\right ]^{m - 2}, \left [ m-1\right ]^1, \left [ m+1\right ]^1 \right\}$ and $\SE(\Gamma_{T_{4m}}) = \frac{4(m-1)(m-3)}{(m+1)} + 4.$	
		\end{enumerate} 
		
\end{cor}

\begin{proof}
	We know that $\frac{T_{4m}}{Z(T_{4m})} \cong D_{2m}$. Therefore, by Theorem \ref{Energy-NCCC2} and  Theorem \ref{Energy-NCCC1}, we get the required result.
\end{proof}

\begin{cor}\label{cor U6m}
Consider the group $U_{6m}$, where $m \geq 2$. Then

 $\Spec(\Gamma_{U_{6m}}) = \left\{ [0]^{2m - 2} ,\left [m\right ]^1, \left [-m\right ]^1\right\}$,
 $\L-Spec(\Gamma_{U_{6m}}) = \left\{ [0]^1, \left [m\right ]^{2m - 2}, \left [ 2m\right ]^1 \right\}$,

 $\Q-Spec(\Gamma_{U_{6m}}) = \left\{ [0]^1, \left [m\right ]^{2m - 2}, \left [ 2m\right ]^1 \right\}$  
and $\E(\Gamma_{U_{6m}}) = \LE(\Gamma_{U_{6m}}) = \SE(\Gamma_{U_{6m}}) = 2m.$
\end{cor}

\begin{proof}
	We know that $\frac{U_{6m}}{Z(U_{6m})} \cong D_{2 \times 3}$. Therefore, by Theorem \ref{Energy-NCCC2}, we get the required result.
\end{proof}

\begin{cor}\label{C3}
Consider the group $U_{(n,m)}$, where $m \geq 3$ and $n \geq 2$.	
	\begin{enumerate}
		\item For $m$ is odd\\
	\rm{(i)} $\Spec(\Gamma_{U_{(n,m)}}) = \left\{ [0]^{\frac{n(m + 1)}{2} - 2} ,\left [n\sqrt{\frac{m - 1}{2}}\right ]^1, \left [-n\sqrt{\frac{m - 1}{2}}\right ]^1\right\}$ and
	$\E(\Gamma_{U_{(n,m)}}) = 2n\sqrt{\frac{m - 1}{2}}.$
	
	\noindent \rm{(ii)} $\L-Spec(\Gamma_{U_{(n,m)}}) = \left\{ [0]^1, \left [n\right ]^{\frac{n(m-1)}{2} - 1}, \left [\frac{n(m - 1)}{2} \right ]^{n - 1}, \left [ \frac{n(m+1)}{2}\right ]^1 \right\}$ and \[\LE(\Gamma_{U_{(n,m)}}) =  \frac{n^2(m-1)(m-3)}{(m+1)} + 2n.\]
	\rm{(iii)} $\Q-Spec(\Gamma_{U_{(n,m)}}) = \left\{ [0]^1, \left [n\right ]^{\frac{n(m-1)}{2} - 1}, \left [\frac{n(m - 1)}{2} \right ]^{n - 1}, \left [ \frac{n(m+1)}{2}\right ]^1 \right\}$ and \[\SE(\Gamma_{U_{(n,m)}}) =  \frac{n^2(m-1)(m-3)}{(m+1)} + 2n.\]	
	\item For $m$ is even and $\frac{m}{2}$ is odd\\
	\rm{(i)} $\Spec(\Gamma_{U_{(n,m)}}) = \left\{ [0]^{\frac{n(m + 2)}{2} - 2} ,\left [n\sqrt{m - 2}\right ]^1, \left [-n\sqrt{m - 2}\right ]^1\right\}$ and
	$\E(\Gamma_{U_{(n,m)}}) = 2n\sqrt{m - 2}.$
	\rm{(ii)} $\L-Spec(\Gamma_{U_{(n,m)}}) = \left\{ [0]^1, \left [2n\right ]^{\frac{n(m-2)}{2} - 1}, \left [\frac{n(m - 2)}{2} \right ]^{2n - 1}, \left [ \frac{n(m+2)}{2}\right ]^1 \right\}$ and \[\LE(\Gamma_{U_{(n,m)}}) =  \frac{2n^2(m-2)(m-6)}{(m+2)} + 4n.\]
	\rm{(iii)} $\Q-Spec(\Gamma_{U_{(n,m)}}) = \left\{ [0]^1, \left [2n\right ]^{\frac{n(m-2)}{2} - 1}, \left [\frac{n(m - 2)}{2} \right ]^{2n - 1}, \left [ \frac{n(m+2)}{2}\right ]^1 \right\}$ and \[\SE(\Gamma_{U_{(n,m)}}) =  \frac{2n^2(m-2)(m-6)}{(m+2)} + 4n.\]
	
	\item For $m$ is even and $\frac{m}{2}$ is even\\	
	\rm{(i)} \small{$\Spec(\Gamma_{U_{(n,m)}}) = \left\{ [0]^{\frac{n(m + 2)}{2} - 3} , \left[ -n \right]^1 ,\left[ \frac{n\left(1 + \sqrt{4m - 7}\right)}{2}\right]^1, \left[ \frac{n\left(1 - \sqrt{4m - 7}\right)}{2}\right]^1 \right\}$ and 
	\[\E(\Gamma_{U_{(n,m)}}) = n\left(1 + \sqrt{4m - 7}\right).\]}
	\rm{(ii)} $\L-Spec(\Gamma_{U_{(n,m)}}) = \left\{ [0]^1, \left [2n\right ]^{\frac{n(m-2)}{2} - 1}, \left [\frac{nm}{2} \right ]^{2n - 2}, \left [ \frac{n(m+2)}{2}\right ]^2 \right\}$ and  \[\LE(\Gamma_{U_{(n,m)}}) =  \frac{2n^2(m-2)(m-4)}{(m+2)} + 4n .\]  
	\rm{(iii)} $\Q-Spec(\Gamma_{U_{(n,m)}}) = \left\{\left [2n\right ]^{\frac{n(m-2)}{2} - 1},  \left [\frac{nm}{2} \right ]^{2n - 2}, \left[\frac{n(m - 2)}{2} \right]^1, \right. \\ \left. \qquad\qquad\qquad\qquad\left[\frac{n\left(m + 6 + \sqrt{(m - 2)(m + 14)}\right)}{4} \right]^1, \left[\frac{n\left(m + 6 - \sqrt{(m - 2)(m + 14)}\right)}{4} \right]^1 \right\}$ and 
	\[\SE(\Gamma_{U_{(n,m)}}) = \begin{cases}
		4n, &\text{for } m = 4\\
		\frac{2n^2(m-2)(m-4)}{(m+2)} + \frac{n(m-2)\left(\sqrt{1+\frac{16}{m-2}}-1\right)}{2}, &\text{for } m \neq 4.
		%\frac{2n^2(m-2)(m-4)}{(m+2)} + \frac{n\left(\sqrt{(m-2)(m +14)}-(m-2)\right)}{2}, &\text{for } m \neq 4.				
	\end{cases}\]	
\end{enumerate}
\end{cor}

\begin{proof}
		We know that $\frac{U_{(n,m)}}{Z(U_{(n,m)})} \cong D_{2 \times \frac{m}{2}}$ or $D_{2m}$ according as $m$ is even or odd. Therefore, by Theorem \ref{Energy-NCCC2} and  Theorem \ref{Energy-NCCC1}, we get the required result.
\end{proof}

\begin{cor}\label{C4}
Consider the group $SD_{8m}$, where $m \geq 2$.
		\begin{enumerate}
		\item For $m$ is even\\
		\rm{(i)} $\Spec(\Gamma_{SD_{8m}}) = \left\{ [0]^{2m - 2} , \left[ -1 \right]^1 ,\left[ \frac{1 + \sqrt{16m - 7}}{2}\right]^1, \left[ \frac{1 - \sqrt{16m - 7}}{2}\right]^1 \right\}$
		%2, 8, 11, 23, 28, 46, 53, 77, 86, 116, 127, 163, 176 upto 200
		 and 
		\[\E(\Gamma_{SD_{8m}}) = 1 + \sqrt{16m - 7}.\]
		\rm{(ii)} $\L-Spec(\Gamma_{SD_{8m}}) = \left\{ [0]^1, \left[2 \right]^{2m - 2}, \left[ 2m + 1 \right]^2 \right\}$ and  $\LE(\Gamma_{SD_{8m}}) =  \frac{4(2m-1)(2m-2)}{(2m+1)} + 4.$
		  
		\noindent\rm{(iii)} $\Q-Spec(\Gamma_{SD_{8m}}) = \left\{\left [2\right ]^{2m - 2}, \left[2m - 1 \right]^1,\left[\frac{\left(2m + 3 + \sqrt{(2m - 1)(2m+7)}\right)}{2} \right]^1,\right. \\ \left.\qquad\qquad\qquad\qquad\qquad\qquad\qquad\qquad\qquad \left[\frac{\left(2m + 3 - \sqrt{(2m - 1)(2m+7)}\right)}{2} \right]^1 \right\}$ and 
		\[\SE(\Gamma_{SD_{8m}}) = \frac{4(2m-1)(2m-2)}{2m+1} +(2m-1)\left(\sqrt{1+\frac{8}{2m-1}}-1\right).\]
	%	\[\SE(\Gamma_{SD_{8m}}) = \frac{4(2m-1)(2m-2)}{2m+1} +\left(\sqrt{(2m-1)(2m + 7)}-(2m-1)\right).\]
		\item For $m$ is odd\\
		\rm{(i)} $\Spec(\Gamma_{SD_{8m}}) = \left\{ [0]^{2m} ,\left [4\sqrt{\frac{m - 1}{2}}\right ]^1, \left [-4\sqrt{\frac{m - 1}{2}}\right ]^1\right\}$ and
		$\E(\Gamma_{SD_{8m}}) = 8\sqrt{\frac{m - 1}{2}}.$
		
		\noindent \rm{(ii)} $\L-Spec(\Gamma_{SD_{8m}}) = \left\{ [0]^1, \left [4\right ]^{2m - 3}, \left [2(m - 1)\right ]^{3}, \left [ 2(m + 1)\right ]^1 \right\}$ and \[\LE(\Gamma_{SD_{8m}}) =  \frac{16(m-1)(m-3)}{(m+1)} + 8.\]
		\rm{(iii)} $\Q-Spec(\Gamma_{SD_{8m}}) = \left\{ [0]^1, \left [4\right ]^{2m - 3}, \left [2(m - 1)\right ]^{3}, \left [ 2(m + 1)\right ]^1 \right\}$  and \[\SE(\Gamma_{SD_{8m}}) =  \frac{16(m-1)(m-3)}{(m+1)} + 8.\]	
	\end{enumerate}

\end{cor}

\begin{proof}
	We know that $\frac{SD_{8m}}{Z(SD_{8m})} \cong D_{2 \times 2m}$ or $D_{2m}$ according as $m$ is even or odd. Therefore, by Theorem \ref{Energy-NCCC2}, we get the required result.
\end{proof}

\begin{thm}\label{Energy-NCCC3}
Consider the group $V_{8m}$, where  $m \geq 2$.
\begin{enumerate}
	\item For $m$ is even\\
	\rm{(i)} $\Spec(\Gamma_{V_{8m}})\! =\! \left\{ [0]^{2m - 1} , \left[ -2 \right]^1 ,\left[ 1 + \sqrt{8m - 7}\right]^1, \left[ 1 - \sqrt{8m - 7}\right]^1 \right\}$ and $$\E(\Gamma_{V_{8m}}) = 2\left( 1 + \sqrt{8m - 7}\right).$$
	\rm{(ii)} $\L-Spec(\Gamma_{V_{8m}}) = \left\{ [0]^1, \left [4\right ]^{2m - 3}, \left [2m \right ]^{2}, \left [2(m+1)\right ]^2 \right\}$ and $\LE(\Gamma_{V_{8m}}) =  \frac{16(m-1)(m-2)}{(m+1)} + 8.$

\noindent	(iii) $\Q-Spec(\Gamma_{V_{8m}})\! =\! \left\{\left [4\right ]^{2m - 3},  \left [2m \right ]^{2}, \left[2(m - 1) \right]^1, \left[m + 3 + \sqrt{(m - 1)(m+7)}\right]^1,\right. \\ \left.
	\qquad\qquad\qquad\qquad\qquad\qquad\qquad\qquad\qquad\left[m + 3 - \sqrt{(m - 1)(m+7)}\right]^1 \right\}$ and $$\SE(\Gamma_{V_{8m}}) = \begin{cases}
			8 &\text{for } m = 2\\
			\frac{16(m-1)(m-2)}{m+1} + 2(m-1)\left(\sqrt{1+\frac{8}{m-1}}-1\right) &\text{for } m \neq 2.
		%		\frac{16(m-1)(m-2)}{m+1} + 2\left(\sqrt{(m-1)(m + 7)}-(m-1)\right) &\text{for } m \neq 2.
		\end{cases}$$
	\item For $m$ is odd\\
		\rm{(i)} $\Spec(\Gamma_{V_{8m}}) = \left\{ [0]^{2m - 2} , \left[ -1 \right]^1 ,\left[\frac{1 + \sqrt{16m - 7}}{2}\right]^1, \left[\frac{1 - \sqrt{16m - 7}}{2}\right]^1 \right\}$ and $$\E(\Gamma_{V_{8m}}) =  1 + \sqrt{16m - 7}.$$
		\rm{(ii)} $\L-Spec(\Gamma_{V_{8m}}) = \left\{ [0]^1, \left [2\right ]^{2m - 2}, \left [2m + 1 \right ]^{2}\right\}$ and $\LE(\Gamma_{V_{8m}}) =  \frac{(4m - 2)^2 + 8}{(2m+1)}.$
		
	\noindent	\rm{(iii)} $\Q-Spec(\Gamma_{V_{8m}}) = \left\{ [2]^{2m - 2},  [2m - 1]^{1}, \left[\frac{2m + 3 + \sqrt{(2m - 1)(2m + 7)}}{2}\right]^1,\right. \\ \left. \qquad\qquad\qquad\qquad\qquad\qquad\qquad\qquad
			\left[\frac{2m + 3 - \sqrt{(2m - 1)(2m + 7)}}{2}\right]^1 \right\}$ and $$\SE(\Gamma_{V_{8m}}) = \frac{3(2m-1)(2m-3)}{2m+1} + (2m-1)\sqrt{1+\frac{8}{(2m-1)}}.$$
		%	$$\SE(\Gamma_{V_{8m}}) = \frac{3(2m-1)(2m-3)}{2m+1} + \sqrt{(2m-1)(2m +7)}.$$
		
	\end{enumerate} 
\end{thm}

\begin{proof}
	\rm{(a)} If $m$ is even then by \cite[Proposition 2.4]{SA-CA-2020},	we have $CCC(V_{8m}) = 2K_2 \cup K_{2m - 2}$. Therefore,
	$\Gamma_{V_{8m}} \cong K_{2 \cdot 2 , 1\cdot 2m - 2}$.
	 Here $|v(\Gamma_{V_{8m}})| = 2(m +1)$ and $r =3$.\\
	\rm{(i)} Using Lemma \ref{Kapbq-poly}(a), we get \begin{align*}
		P(\Gamma_{V_{8m}}, x)
		%&= x^{2(m +1) - 3}(x+2)^{2-1}(x+(2m - 2))^{1-1}\\
		%&\left(x^2+ (2(1 - 2) + (2m - 2)(1 - 1))x + 2(2m - 2)(1 - 2 - 1)\right).\\
		&= x^{2(m + 1)- 3}\left(x + 2\right) \left(x^2- 2x - 8(m - 1)\right).
	\end{align*}
	Thus, $\Spec(\Gamma_{V_{8m}}) = \left\{ [0]^{2m - 1} , \left[ -2 \right]^1 ,\left[ 1 + \sqrt{1+8(m - 1)}\right]^1, \left[ 1 - \sqrt{1+8(m - 1)}\right]^1 \right\}$ and so
	\small{
	\begin{align*}
		\E(\Gamma_{V_{8m}}) &= (2m - 1) \times 0 + 1 \times \left|-2\right| + 1 \times \left| 1 + \sqrt{1+8(m - 1)} \right| + 1 \times \left| 1 - \sqrt{1+8(m - 1)} \right|\\
		%&= 2 + \left( 1 + \sqrt{1+8(m - 1)}\right) + \left(\sqrt{1+8(m - 1)} - 1\right)\\
		&= 2\left( 1 + \sqrt{1+8(m - 1)}\right).
	\end{align*}
}
	\rm{(ii)} Using Lemma \ref{Kapbq-poly}(b), we get \begin{align*}
		P_L(\Gamma_{V_{8m}}, x) 
		%&= x(x - (2\times2 + (2m - 2)(1 - 1)))^{(2m - 2-1)}\\
		%&(x - (2m - 2 + 2(2 - 1)))^{2(2-1)} (x - 2(m +1))^{2 + 1 - 1}\\
			&= x(x - 4)^{2m - 3} \left(x - 2m \right)^{2} \left(x - 2(m+1) \right)^2 .
	\end{align*}
	
	So, $\L-Spec(\Gamma_{V_{8m}}) = \left\{ [0]^1, \left [4\right ]^{2m - 3}, \left [2m \right ]^{2}, \left [2(m+1)\right ]^2 \right\}$ .
	
	Now $|2e(\Gamma_{V_{8m}})|=  8\left (2m - 1\right )$ and so $\Delta(\Gamma_{V_{8m}}) = \frac{|2e(\Gamma_{V_{8m}})|}{|v(\Gamma_{V_{8m}})|} = \frac{4(2m-1)}{(m+1)}$. We have
	
	\[
	\left|0 - \Delta(\Gamma_{V_{8m}}) \right| =  \left|  \frac{-4(2m-1)}{(m+1)} \right| =  \frac{4(2m-1)}{(m+1)}, \quad
		\left|4 - \Delta(\Gamma_{V_{8m}}) \right| =  \left|\frac{-4(m - 2)}{(m+1)}\right| =  \frac{4(m - 2)}{(m+1)},
\]
	
	\begin{align*}
		\left|2m - \Delta(\Gamma_{V_{8m}})\right| &=  \left| \frac{2(m - 2)(m-1)}{(m+1)}\right| =  \frac{2(m - 2)(m-1)}{(m+1)}
	\end{align*}
	and
	\begin{align*}
		\left|(2m + 2) - \Delta(\Gamma_{V_{8m}}) \right| &=  \left|  \frac{2\left((m-1)^2+2\right)}{(m+1)}\right| =  \frac{2((m-1)^2+2)}{(m+1)}.
	\end{align*}
	Therefore,
	\begin{align*}
		\LE(\Gamma_{V_{8m}}) &= 1 \times \left|  \frac{-4(2m-1)}{(m+1)} \right| + \left(2m - 3\right) \times \left|\frac{-4(m - 2)}{(m+1)}\right|
		+ 2 \times \left| \frac{2(m - 2)(m-1)}{(m+1)}\right| + 2 \times \left|  \frac{2\left((m-1)^2+2\right)}{(m+1)}\right|\\
		&=  \frac{16(m-1)(m-2)}{(m+1)} + 8 .
	\end{align*}
	
\noindent	\rm{(iii)} Using Lemma \ref{Kapbq-poly}(c), we get
		\begin{align*}
		P_Q(\Gamma_{V_{8m}}, x) 
		%= \prod_{i=1}^{2}\left(x - \left( 2m + 2 - p_i\right)\right)^{a_i(p_i-1)} \prod_{i=1}^{2}\left(x - \left( 2m + 2 - 2p_i\right)\right)^{a_i - 1}\\
	%	&\left(x^2- x\left(2(2m + 2)- p_2\right)+ \left(2m + 2\right)^2 - 2(2m + 2)(p_1 + p_2) \right. \\ & \left. \qquad\qquad\qquad\qquad\qquad\qquad\qquad\qquad + (2m + 2)(2p_1 + p_2) - 2p_1p_1\right)\\
		= &\left(x - 4\right)^{2m - 3} \left(x - 2m \right)^{2} \left(x - 2(m - 1) \right)
		 \left(x^2- \left(2m + 6\right)x 
		  + 4m^2 + 4m + 16\right)
		\end{align*}
	Therefore  
	\begin{align*}
		\Q-Spec(\Gamma_{V_{8m}}) = \biggl\{\left [4\right ]^{2m - 3},  \left [2m \right ]^{2}, \left[2(m - 1) \right]^1, &\left[m + 3 + \sqrt{(m - 1)(m+7)}\right]^1,
		\left[m + 3 - \sqrt{(m - 1)(m+7)}\right]^1 \biggr\}.
	\end{align*}
		
	We have
	
	\begin{align*}
		\left|2(m-1) - \Delta(\Gamma_{V_{8m}}) \right| &=  \left| \frac{2\left((m-2)^2-3\right)}{(m+1)}\right| = \begin{cases}
			2 &\text{for } m = 2\\
		    \frac{2\left((m-2)^2-3\right)}{(m+1)} &\text{for } m \geq 4,
		\end{cases}
	\end{align*}
\begin{align*}
	\left|\left(m + 3 + (m - 1)\sqrt{1 + \frac{8}{m - 1}}\right) \!-\! \Delta(\Gamma_{V_{8m}})  \right| &\!\!=\!\!  \left| \frac{\left((m - 2)^2 + 3 + (m + 1)\sqrt{(m - 1)(m+7)}\right)}{(m+1)} \right|\\ 
%	\!\!&=\!\! \begin{cases}
%		4 &\text{for } m = 2\\
	   &=\frac{\left((m - 2)^2 + 3 + (m + 1)\sqrt{(m - 1)(m+7)}\right)}{(m+1)} 
	   %&\text{for } m \geq 4.
%
\end{align*}
and
\begin{align*}
	\left|\left(m + 3 - (m - 1)\sqrt{1 + \frac{8}{m - 1}}\right)\! -\! \Delta(\Gamma_{V_{8m}}) \right|\!\! &=\!\!  \left|\frac{\left((m - 2)^2 + 3 - (m + 1)\sqrt{(m - 1)(m+7)}\right)}{(m+1)} \right|.
\\ 
%	&\!\!=\!\! \begin{cases}
%		2 &\text{for } m = 2\\
%		\frac{-\left((m - 2)^2 + 3 - (m + 1)\sqrt{(m - 1)(m+7)}\right)}{(m+1)}  &\text{for } m \geq 4,
%	\end{cases}
\end{align*}
Let $f_1(m) := (m - 2)^2 + 3 - (m + 1)\sqrt{(m - 1)(m+7)}$. Then for $m \geq 4$, $((m - 2)^2 + 3)^2 - (m + 1)^2(m - 1)(m+7) = -8m^2(2m - 3) - 48m + 56 < 0$ and so $(m - 2)^2 + 3 < (m + 1)\sqrt{(m - 1)(m+7)} \implies f_1(m) < 0$. Also we have, $f_1(2) = - 6$. Therefore,
\small{
\begin{align*}
	\left|\left(m + 3 - (m - 1)\sqrt{1 + \frac{8}{m - 1}}\right)\! -\! \Delta(\Gamma_{V_{8m}}) \right| 
	\!\!=\!\! \frac{-\left((m - 2)^2 + 3 - (m + 1)\sqrt{(m - 1)(m+7)}\right)}{(m+1)}.
\end{align*}
}
%since for $m \geq 4$, $((m - 2)^2 + 3)^2 - (m + 1)^2(m - 1)(m+7) = -8m^2(2m - 3) - 48m + 56 < 0$ and so $(m - 2)^2 + 3 < (m + 1)\sqrt{(m - 1)(m+7)} < 0$. Therefore
	\begin{align*}
		\SE(\Gamma_{V_{8m}}) &= \left( 2m - 3 \right) \times \left|\frac{-4(m - 2)}{(m+1)}\right| + 2 \times \left| \frac{2(m - 2)(m-1)}{(m+1)}\right|\\
		&+ 1 \times \left| \frac{2\left((m-2)^2-3\right)}{(m+1)}\right|
		+\left| \frac{\left((m - 2)^2 + 3 + (m + 1)\sqrt{(m - 1)(m+7)}\right)}{(m+1)} \right|\\
		&\qquad\qquad\qquad\qquad\qquad\quad + \left|\frac{\left((m - 2)^2 + 3 - (m + 1)\sqrt{(m - 1)(m+7)}\right)}{(m+1)} \right|\\
		&= \begin{cases}
			8 &\text{for } m = 2\\
			\frac{16(m-1)(m-2)}{(m+1)} + 2(m-1)\left(\sqrt{1+\frac{8}{m-1}}-1\right) &\text{for } m \geq 4.
				%		\frac{16(m-1)(m-2)}{m+1} + 2\left(\sqrt{(m-1)(m + 7)}-(m-1)\right) &\text{for } m \neq 2. 
		\end{cases} 
	\end{align*}
	Hence the result follows.
	
\noindent	\rm{(b)} If $m$ is odd, then by \cite[Proposition 2.4]{SA-CA-2020}, we have $\Gamma_{V_{8m}} \cong K_{2 \cdot 1 , 2m - 1} \cong \Gamma_{D_{2\times 4m}}$. Hence, the result follows from Corollary \ref{C1}.
\end{proof}

\section{Groups whose NCCC-graphs are integral, L-integral and Q-integral}
In this section, we determine whether NCCC-graphs of the groups considered in  Section 2 are integral, L-integral and Q-integral. In view of Theorem \ref{Energy-NCCC1}--Theorem \ref{Energy-NCCC3}, we get the following results.
\begin{thm}
	Let $G$ be a finite non-abelian group with center $Z(G)$.
\begin{enumerate}
\item Suppose that $G = D_{2m}$, $T_{4m}$, $U_{(n,m)}$, $SD_{8m}$, where $m$ is odd. Then $\Gamma_{G}$ is integral if and only if $m = 2a^2 + 1$, where $a \geq 1$.

\item Suppose that $G = D_{2m}$, $U_{(n,m)}$, where $m$ is even and $\frac{m}{2}$ is odd. Then $\Gamma_{G}$ is integral if and only if $m = a^2 + 2$, where $a \geq 2$.

\item  Suppose that $G = D_{2m}$, where $m$  and  $\frac{m}{2}$ are even. Then $\Gamma_{G}$ is integral if and only if $4m - 7$ is a perfect square.

\item  Suppose that $m$  is  an even integer. Then  $\Gamma_{T_{4m}}$ and $\Gamma_{SD_{8m}}$ are integral if and only if $8m -7$ and $16m -7$ are perfect squares respectively.

%\item  If $m$ is even and $8m -7$ is a perfect square then $\Gamma_{V_{8m}}$ is integral. If $m$ is odd and $16m -7$ is a perfect square then $\Gamma_{V_{8m}}$ is integral.

\item  For $m$ is even, $\Gamma_{V_{8m}}$ is integral if and only if $8m -7$ is a perfect square. For $m$ is odd, $\Gamma_{V_{8m}}$ is integral if and only if $16m -7$ is a perfect square.
\end{enumerate}   
\end{thm}
\begin{thm}
Let $G$ be a finite non-abelian group with center $Z(G)$. Then the NCCC-graph of $G$ is
\begin{enumerate}
\item  integral, L-integral and Q-integral if $\frac{G}{Z(G)} \cong \mathbb{Z}_p \times \mathbb{Z}_p$ or $G \cong U_{6m}$.
\item L-integral if $G \cong D_{2m}$, $T_{4m}$, $SD_{8m}$, $U_{(n,m)}$ or $V_{8m}$. In general, $\Gamma_{G}$ is  L-integral if  $\frac{G}{Z(G)} \cong D_{2m}$.
\item Q-integral if $G \cong D_{2m}$, $T_{4m}$, $SD_{8m}$ or $U_{(n,m)}$, where $m$ is odd.
\item  Q-integral if $G \cong D_{2m}$ or $U_{(n,m)}$, where $m$ is even and $\frac{m}{2}$ is odd.
\end{enumerate} 
\end{thm}
The following lemma is useful for the cases when $m$  and $\frac{m}{2}$ are both even.
%\begin{thm}
%	Let $G$ be a finite non-abelian group with centre $Z(G)$. Then the NCCC graph of $G$ is 
%	\begin{enumerate}
%		\item integral if $\frac{G}{Z(G)} \cong \mathbb{Z}_p \times \mathbb{Z}_p.$
%		\item L- integral if $\frac{G}{Z(G)} \cong \mathbb{Z}_p \times \mathbb{Z}_p,$ $G \cong D_{2m}$, $T_{4m}$, $SD_{8m}$, $U_{6m}$,  $U{(n,m)}$ and $V_{8m}$.
%		\item Q- integral if $\frac{G}{Z(G)} \cong \mathbb{Z}_p \times \mathbb{Z}_p.$
%	\end{enumerate}
%\end{thm}
\begin{lem}\label{L1}
	Let $m$ be any positive integer. Then
	\begin{enumerate}
		\item $m^2 + 6m - 7$ is a perfect square if and only if $m = 1, 2.$
		\item $m^2 + 12m - 28$ is a perfect square if and only if $m = 2, 4, 11.$
		\item $4m^2 + 12m - 7$ is a perfect square if and only if $m = 1.$
	\end{enumerate}
\end{lem}

\begin{proof}
	
	(a) Let $m^2 + 6m - 7$ be a perfect square. Then there exist integers $k$ such that $m^2 + 6m - 7 = k^2$ which in turn gives $(m + 3 + k)(m + 3 - k) = 16.$ Therefore, we have the following cases.
	
	\noindent\textbf{Case 1.} $m + 3 + k = 16$ and $m + 3 - k = 1$\\
	In this case, we have $m + k =13$ and $m - k = -2$ which gives $m = \frac{11}{2}$ and $k = \frac{15}{2}$.
	
	\noindent\textbf{Case 2.} $m + 3 + k = -16$ and $m + 3 - k = -1$\\
	In this case, we have $m + k = -19$ and $m - k = -4$ which gives $m = \frac{-23}{2}$ and $k = \frac{-15}{2}$.
	
	\noindent\textbf{Case 3.} $m + 3 + k = 1$ and $m + 3 - k = 16$\\
	In this case, we have $m + k = -2$ and $m - k = 13$ which gives $m = \frac{11}{2}$ and $k = \frac{-15}{2}$.
	
\noindent	\textbf{Case 4.} $m + 3 + k = -1$ and $m + 3 - k = -16$\\
	In this case, we have $m + k = -4$ and $m - k = -19$ which gives $m = \frac{-23}{2}$ and $k = \frac{15}{2}$.
	
\noindent\textbf{Case 5.} $m + 3 + k = 2$ and $m + 3 - k = 8$\\
	In this case, we have $m + k = -1$ and $m - k = 5$ which gives $m = 2$ and $k = -3$.
	
\noindent	\textbf{Case 6.} $m + 3 + k = -2$ and $m + 3 - k = -8$\\
	In this case, we have $m + k = -5$ and $m - k = -11$ which gives $m = -8$ and $k = 3$.
	
\noindent	\textbf{Case 7.} $m + 3 + k = 8$ and $m + 3 - k = 2$\\
	In this case, we have $m + k = 5$ and $m - k = -1$ which gives $m = 2$ and $k = 3$.\\
	\textbf{Case 8.} $m + 3 + k = -8$ and $m + 3 - k = -2$\\
	In this case, we have $m + k = -11$ and $m - k = -5$ which gives $m = -8$ and $k = -3$.
	
\noindent	\textbf{Case 9.} $m + 3 + k = 4$ and $m + 3 - k = 4$\\
	In this case, we have $m + k = 1$ and $m - k = 1$ which gives $m = 1$ and $k = 0$.
	
\noindent	\textbf{Case 10.} $m + 3 + k = -4$ and $m + 3 - k = -4$\\
	In this case, we have $m + k = -7$ and $m - k = -7$ which gives $m = -7$ and $k = 0$.\\
	Hence the result follows.
	
	(b) If $m^2 + 12m - 28$ is a perfect square. Then there exist integers $k$ such that $m^2 + 12m - 28 = k^2$ which in turn gives $(m + 6 + k)(m + 6 - k) = 64.$ By considering different cases as above we get $m = 2, 4, 11$. Hence, the result follows.
	
	(c) If $4m^2 + 12m - 7$ is a perfect square. Then there exist integers $k$ such that $4m^2 + 12m - 7 = k^2$ which in turn gives $(2m + 3 + k)(2m + 3 - k) = 16.$ By considering different cases as above we get $m = 1$. Hence, the result follows.
\end{proof}
\begin{thm}
Let $G$ be a finite non-abelian group.
\begin{enumerate}
\item  If $G \cong D_{2m}$ or $U_{(n,m)}$,  where  $m$ and 	$\frac{m}{2}$ are both even, then $\Gamma_{G}$ is Q-integral if and only if $m = 4$.
\item If $G \cong T_{4m}$, where  $m$ is even, then $\Gamma_{G}$ is Q-integral if and only if $m = 2$.
\item If $G \cong V_{8m}$ then $\Gamma_{G}$ is Q-integral if and only if $m = 2$.
\item If $G \cong SD_{8m}$, where $m$ is even, then $\Gamma_{G}$ is not Q-integral.
\end{enumerate}
\end{thm}
\begin{proof}
(a) Suppose that $m$ and $\frac{m}{2}$ are both even. If $G \cong D_{2m}$ (where $m \geq 4$) then, in view of Corollary \ref{C1}(c), it is sufficient to show that $\frac{m + 6 + \sqrt{m^2 + 12m - 28}}{4}$ and $\frac{m + 6 - \sqrt{(m - 2)(m + 14)}}{4}$ are integers. By Lemma \ref{L1}(b), we have $m = 4$. Note that   $\frac{m + 6 + \sqrt{m^2 + 12m - 28}}{4} = 4$ and $\frac{m + 6 - \sqrt{m^2 + 12m - 28}}{4} = 1$. 

 If $G \cong U_{(n,m)}$ (where $m \geq 4$ and $n \geq 2$) then, in view of Corollary \ref{C3}(c), it is sufficient to show that $\frac{n\left(m + 6 + \sqrt{m^2 + 12m  - 28}\right)}{4}$ and $\frac{n\left(m + 6 - \sqrt{m^2 + 12m  - 28}\right)}{4}$ are integers. By Lemma \ref{L1}(b), we have $m = 4$. Note that 
	  $\frac{n\left(m + 6 + \sqrt{m^2 + 12m  - 28}\right)}{4} = 4n$ and $\frac{n\left(m + 6 - \sqrt{m^2 + 12m  - 28}\right)}{4} = n$.

(b) Suppose that $m$ is even. If $G \cong T_{4m}$ (where $m \geq 2)$ then, in view of Corollary \ref{C2}(a), it is sufficient to show that $\frac{m + 3 + \sqrt{m^2 + 6m -7}}{2}$ and $\frac{m + 3 - \sqrt{m^2 + 6m -7}}{2}$ are integers. By Lemma \ref{L1}(a), we have $m = 2$. Note that  $\frac{m + 3 + \sqrt{m^2 + 6m -7}}{2}$ $= 4$ and $\frac{m + 3 - \sqrt{m^2 + 6m -7}}{2} = 1$. 

(c) Suppose that $m$ is even. If $G \cong V_{8m}$ (where $m \geq 2)$ then, in view of  Theorem \ref{Energy-NCCC3}(a), it is sufficient to show that $m + 3 + \sqrt{m^2 + 6m -7}$ and $m + 3 - \sqrt{m^2 + 6m -7}$ are integers. By Lemma \ref{L1}(a), we have $m = 2$. Note that  $m + 3 + \sqrt{m^2 + 6m -7} = 8$ and $m + 3 - \sqrt{m^2 + 6m -7} = 2$. 

If $m$ is odd then, by Theorem \ref{Energy-NCCC3}(b), it is sufficient to show that $\frac{2m + 3 + \sqrt{4m^2 + 12m - 7}}{2}$ and $\frac{2m + 3 - \sqrt{4m^2 + 12m - 7}}{2}$ are integers. By Lemma \ref{L1}(c), we have $m = 1$ but $m \geq 3$. Thus, $\Gamma_{G}$ is not Q-integral.

(d) If $m$ is even and $G = SD_{8m}$ (where $m\geq 2$)
 then, by Corollary \ref{C4}(a), it is sufficient to show that $\frac{2m + 3 + \sqrt{4m^2 + 12m - 7}}{2}$ and $\frac{2m + 3 - \sqrt{4m^2 + 12m - 7}}{2}$ are integers. By Lemma \ref{L1}(c), we have $m = 1$. Thus, $\Gamma_{G}$ is not Q-integral.
\end{proof}
In view of Theorem \ref{Energy-NCCC2} and Lemma \ref{L1}(a), we also have the following result.
\begin{thm}
Let $G$ be a finite group and $\frac{G}{Z(G)} \cong D_{2m}$, where $m \geq 3$. Then $\Gamma_{G}$ is  Q-integral if and only if $m$ is odd.
%\begin{enumerate}
%\item $\Gamma_{G}$ is  Q-integral if $m$ is odd.
%\item If $m$ is even then $\Gamma_{G}$ is  Q-integral if and only if $m =2$ and $|Z(G)|$ is even.
%\end{enumerate}  
\end{thm}

\section{Groups whose NCCC-graphs are hyperenergetic, L-hyperenergetic and Q-hyperenergetic}
In this section, we consider  NCCC-graphs of the groups  considered in Section 2 and determine whether those graphs are hyperenergetic, L-hyperenergetic and Q-hyperenergetic. 
\begin{thm}
	Let $G$ be a finite non-abelian group with center $Z(G)$ such that $\frac{G}{Z(G)} \cong \mathbb{Z}_p \times \mathbb{Z}_p$, where $p$ is any prime. Then the NCCC-graph of $G$ is neither hyperenergetic, L-hyperenergetic nor Q-hyperenergetic.
\end{thm}

\begin{proof}
	Let $|Z(G)| = z$ and $n = \frac{(p-1)z}{p}$. Then $|v(\Gamma_G)|=n(p+1)$ and, by Theorem \ref{Energy-NCCC1} and \eqref{Kn}, we get
	\[
	\E(\Gamma_G) - \E(K_{|v(\Gamma_G)|}) = -2(n - 1). 
	\]
	Hence, $\E(\Gamma_G) - \E(K_{|v(\Gamma_G)|}) \leq 0$. By Theorem \ref{Energy-NCCC1}, we also have $\E(\Gamma_G) = \LE(\Gamma_G) = \SE(\Gamma_G)$ and so the NCCC-graph of $G$ is neither hyperenergetic, L-hyperenergetic nor Q-hyperenergetic.
\end{proof}
The following corollary follows immediately.
\begin{cor}
	Let $G$ be a non-abelian group of order $p^m$ with $|Z(G)| = p^{m-2}$. Then the NCCC-graph of $G$ is neither hyperenergetic, L-hyperenergetic nor Q-hyperenergetic.
\end{cor}

\begin{thm}
	Let $G$ be a finite group with center $Z(G)$ and $\left|Z(G)\right|$ = $z$. If $\frac{G}{Z(G)}$ is isomorphic to   $D_{2m}$ (where $m \geq 3$). Then the NCCC-graph of $G$ is
	\begin{enumerate} 
		\item borderenergetic if $m = 3$ and $z = 1$.
		\item not hyperenergetic.
		\item L-borderenergetic if $m = 3$ and $z = 1$.
		\item L-hyperenergetic except for $m \geq 5$ and $z = 1$ ; $m = 3$ and $z \geq 2$; $m = 5$ and $z = 1,2$ and $m = 7$ and $z = 1$.  
		\item Q-borderenergetic if $m = 3$ and $z = 1$.
		\item Q-hyperenergetic except for $m = 4$ and $z = 2$; $m \geq 5$ and $z = 1$ ; $m = 3$ and $z \geq 2$; $m = 5$ and $z = 1,2$ and $m = 7$ and $z = 1$.
	\end{enumerate}
\end{thm}

\begin{proof}
We have $|v(\Gamma_G)| = \frac{z(m +1)}{2}$. Consider the following two cases.

\noindent\textbf{Case 1.} $m$ is even\\
In this case, $m \geq 4$ and $z \geq 2$. By Theorem \ref{Energy-NCCC2} and \eqref{Kn}, we have
	\[
	\E(\Gamma_G) - \E(K_{|v(\Gamma_G)|}) = \frac{z\left(-1 + \sqrt{8m-7} - 2m\right) + 4}{2}< 0, 
	\]
	since $\left(-1 + \sqrt{8m-7} - 2m\right) \leq -4$ for $m \geq 4$.
	
	By Theorem \ref{Energy-NCCC2} and \eqref{Kn}, we have 
	\[
	\LE(\Gamma_G) - \LE(K_{|v(\Gamma_G)|}) = \frac{f_1(m,z)}{m+1}, 
	\]
where $f_1(m,z) = m^2 z^2-m^2 z-3 m z^2+2 m+2 z^2+z+2.$ Then $f_1(m,z) = \frac{zm^2(z - 2)}{2} + \frac{z^2m(m - 6)}{2} + 2 m+2 z^2+z+2 > 0$ for $m \geq 6$ and $z \geq 2.$ Now, $f_1(4,z) = 6z^2 - 15z + 10 > 0 $ for $z \geq 2.$
	Hence, $$\LE(\Gamma_G) - \LE(K_{|v(\Gamma_G)|}) > 0.$$
	
	By Theorem \ref{Energy-NCCC2} and \eqref{Kn}, we have
	\begin{align*}
		\SE(\Gamma_G) - \SE(K_{|v(\Gamma_G)|})= \frac{f_2(m,z)}{2(m+1)},
	\end{align*} 
%	$\SE(\Gamma_G) - \SE(K_{|v(\Gamma_G)|})$ 
%	\[
%	= \frac{2z^2m^2 -6z^2m-3zm^2 - z - 4mz +4 z^2+4 m+4 + z(m +1)\sqrt{(m-1)(m+7)}}{2(m+1)}. 
%	\]
where $f_2(m,z) = 2z^2m^2 -6z^2m-3zm^2 - z - 4mz +4 z^2+4 m+4 + z(m +1)\sqrt{(m-1)(m+7)}$. Then $f_2(m,z) = z^2m(m - 6) + zm^2(z - 4) + zm(m - 4) + \frac{z(z - 4)}{4}+4 m+4 + z(m +1)\sqrt{(m-1)(m+7)} > 0$ for $m \geq 6$ and $z \geq 4.$
	
	We have $f_2(4,z) = 12z^2 + (5\sqrt{33} - 65)z + 20 > 0$ for $z \geq 3.$ Also $f_2(4,2) = -4.55.$ Therefore  $f_2(4,z) < 0$  or $> 0$ according as $z =2$ or $z \geq 3$.
	
	We have $f_2(m,2) = 2m^2 - 28m + 18 + 2(m+1)\sqrt{(m-1)(m+7)} < 0$  or $> 0$ according as $m =4$ or $m \geq 6$. Also, $f_2(m,3) = 9m^2 - 62m + 37 + 3(m+1)\sqrt{(m-1)(m+7)} > 0$ for $m \geq 4.$	Hence, 
	\[\SE(\Gamma_G) - \SE(K_{|v(\Gamma_G)|})  \begin{cases}
		< 0, &\text{if }m = 4 \text{ and } z = 2\\
		> 0, &\text{otherwise.}
	\end{cases}\]
%	$\SE(\Gamma_G) - \SE(K_{|v(\Gamma_G)|}) < 0$ if $m = 4$ \& $z = 2.$ Otherwise, $\SE(\Gamma_G) - \SE(K_{|v(\Gamma_G)|}) > 0.$\\
\noindent\textbf{Case 2.  } $m$ is odd
	
	By Theorem \ref{Energy-NCCC2} and \eqref{Kn}, we have
	\[
	\E(\Gamma_G) - \E(K_{|v(\Gamma_G)|}) = z\left(\sqrt{2(m-1)} - m - 1\right) + 2. 
	\]
	We have $\left(\sqrt{2(m-1)} - m - 1\right) < -2$ for $m \geq 5.$
	Hence, $\E(\Gamma_G) - \E(K_{|v(\Gamma_G)|}) = 0$ for $m = 3$ and $z = 1$; otherwise $\E(\Gamma_G) - \E(K_{|v(\Gamma_G)|}) < 0$.
	
	By Theorem \ref{Energy-NCCC2} and \eqref{Kn}, we have 
	\[
	\LE(\Gamma_G) - \LE(K_{|v(\Gamma_G)|}) = \frac{f_3(m,z)}{m+1}, 
	\]
	where $f_3(m,z) = m^2 z^2-m^2 z-4 m z^2+2 m+3 z^2+z+2.$ Then $f_3(m,z) = \frac{zm^2(z - 2)}{2} + \frac{z^2m(m - 8)}{2} + 2 m+3 z^2+z+2 > 0$ for $m \geq 8$ and $z \geq 2$. Now, $f_3(m,1) = -2m + 6 < 0$ for $m \geq 5.$ Therefore, $f_3(m,1) = 0$ or $< 0$ according as $m = 3$ or $m \geq 5$.
	We have $f_3(3,z) = -8z + 8 < 0$ for $z \geq 2.$ Also $f_3(5,z) = 8z^2 - 24z + 12 > 0$ for $z \geq 3$ and $f_3(5,1) = f_3(5,2) = -4.$ So, $f_3(5,z) > 0$ or $< 0$ according as $z \geq 3$ or $z = 1,2$. Again $f_3(7,z) = 24z^2 - 48z + 16 > 0$ for $z \geq 2$ and $f_3(7,1) = -8$. Therefore, $f_3(7,z) > 0$ or $< 0$ according as $z \geq 2$ or $z = 1$.
	Hence, 
	\[\LE(\Gamma_G) - \LE(K_{|v(\Gamma_G)|}) \begin{cases}
		= 0, &\text{if }m = 3 \& z = 1\\
		< 0, &\text{if }m \geq 5 \,\&\, z = 1;~m = 3\,\&\,z \geq 2;\\
		&m = 5 \,\&\, z = 1,2 \text{ and } m = 7 \,\&\,z = 1\\
		> 0, &\text{otherwise.}
	\end{cases}\]
%	$\LE(\Gamma_G) - \LE(K_{|v(\Gamma_G)|}) = 0$ if  $m = 3$ \& $z = 1$, $\LE(\Gamma_G) - \LE(K_{|v(\Gamma_G)|}) < 0$ if $m \geq 5$ \& $z = 1$; $m = 3$ \& $z \geq 2$; $m = 5$ \& $z = 1,2$ and $m = 7$ \& $z = 1$. Otherwise, $\LE(\Gamma_G) - \LE(K_{|v(\Gamma_G)|}) > 0.$\\
	%	Hence, 
	%\[
	%\LE(\Gamma_G) - \LE(K_{|v(\Gamma_G)|})
	%\] = \begin{cases}
		%=0 &\text{for } $m = 3$ \& $z = 1$\\
		%<0 &\text{for } $m \geq 5$ \& $z = 1$; $m = 3$ \& $z \geq 2$; $m = 5$ \& $z %= 1,2$; $m = 7$ \& $z = 1$\\
		%		>0 &\text{otherwise }.
		%	\end{cases}
	By Theorem \ref{Energy-NCCC2}, we also have $\LE(\Gamma_G) = \SE(\Gamma_G)$. Hence the result follows.	
\end{proof}
%\begin{cor}
%	The NCCC graph of $U_{6m} (m \geq 2)$ is neither hyperenergetic, L-hyperenergetic nor Q-hyperenergetic.
%\end{cor}
As a corollary of the above theorem we have the following result. 
\begin{cor}
	Let $G = D_{2m}, T_{4m}, SD_{8m}$ and $U_{(n,m)}$. Then 
	\begin{enumerate}
		\item $\Gamma_{G}$ is borderenergetic if and only  if $G = D_6, D_8, T_8.$
		\item $\Gamma_{G}$ is  not hyperenergetic.
		\item $\Gamma_{G}$ is  L-borderenergetic if and only  if $G = D_6, D_8, T_8.$ 
		\item $\Gamma_{G}$ is not L-hyperenergetic if and only  if $G = D_{2m}$ for $m$ is odd and $m \geq 5$, $D_{12}$, $D_{20}$, $T_{12}, T_{20}$, $SD_{24}$, $U_{(n,m)}$ for $m= 3,4,6$ and $n \geq 2$ and $m=5$ and $n= 2$. 
		\item $\Gamma_{G}$ is  Q-borderenergetic if and only  if $G = D_6$, $D_8$, $T_8$.
		\item $\Gamma_{G}$ is not Q-hyperenergetic if and only  if $G = D_{2m}$ for $m$ is odd and $m \geq 5$, $D_{12}, D_{16}, D_{20}$,  $T_{12}, T_{16}, T_{20}$, $SD_{16}, SD_{24}$, $U_{(n,m)}$ for $m= 3,4,6$ and $n \geq 2$ and $m=5$ and $n= 2$. 
	\end{enumerate}
\end{cor}
We conclude this section with the following result.
\begin{thm}
	The NCCC-graph of $V_{8m} (m \geq 2)$ is L-hyperenergetic and Q- hyperenergetic except for $m = 2$ but not hyperenergetic.
\end{thm}
\begin{proof}
If $m = 2$ then we have $|v(\Gamma_G)| = 6.$ By Theorem \ref{Energy-NCCC3} and \eqref{Kn}, we get
		\[
	\E(\Gamma_G) - \E(K_{|v(\Gamma_G)|}) = -2 < 0. 
	\]
Also, $\E(\Gamma_G) = \LE(\Gamma_G) = \SE(\Gamma_G)$. Thus, $\Gamma_G$ is neither  hyperenergetic,  L- hyperenergetic nor Q-hyperenergetic.
	
If $m$ is even and $m \geq 4$ then we have $|v(\Gamma_G)| = 2m + 2.$ By Theorem \ref{Energy-NCCC3} and (\ref{Kn}), we get
	\[
	\E(\Gamma_G) - \E(K_{|v(\Gamma_G)|}) = 2\left(\sqrt{8m-7} - 2m\right) < 0,
	\]
	\[
	\LE(\Gamma_G) - \LE(K_{|v(\Gamma_G)|}) = \frac{12m^2 - 46m + 38}{m+1} > 0 
	\]
	and
	\[
	\SE(\Gamma_G) - \SE(K_{|v(\Gamma_G)|}) 	= \frac{10m^2 - 54m + 32 + 2(m + 1)\sqrt{(m-1)(m+7)}}{m+1} > 0. 
	\]
Thus, $\Gamma_G$ is not  hyperenergetic but  L-hyperenergetic and Q-hyperenergetic.

If $m$ is odd then we have $|v(\Gamma_G)| = 2m + 1.$ By Theorem \ref{Energy-NCCC3} and \eqref{Kn}, we get
	\[
	\E(\Gamma_G) - \E(K_{|v(\Gamma_G)|}) = 1+ \sqrt{16m-7} - 4m < 0, 
	\]
	\[
	\LE(\Gamma_G) - \LE(K_{|v(\Gamma_G)|}) = \frac{8m^2 - 20m + 12}{2m+1} > 0 
	\]
	and
	\[
	\SE(\Gamma_G) - \SE(K_{|v(\Gamma_G)|}) 	= \frac{4m^2 - 28m + 9 + (2m + 1)\sqrt{(2m-1)(2m+7)}}{2m+1} > 0. 
	\]
Thus, $\Gamma_G$ is not  hyperenergetic but  L-hyperenergetic and Q-hyperenergetic.
\end{proof}

\section{Comparing different energies}
Comparison among various energies become popular because of the E-LE conjecture: 	$\E(\mathcal{G}) \leq \LE(\mathcal{G})$ of Gutman et al. \cite{GAVBR} which is eventually proved to be false in \cite{Liu-09,SSM-09}.  %conjectured that
%\begin{conj}\label{conj}
%	$\E(\mathcal{G}) \leq \LE(\mathcal{G})$.
%\end{conj}
%However it was disproved in \cite{Liu-09,SSM-09}.
% But from then people become interested in comparing different energies. 
Das et al. \cite{DaAo-2018} posed the following  problem comparing $\LE(\mathcal{G})$ and $\SE(\mathcal{G})$.
\begin{prob}\label{Prob 1}\cite[Problem 1]{DaAo-2018}
	Characterize all the graphs for which $\LE(\mathcal{G}) > \SE(\mathcal{G})$, $\LE(\mathcal{G}) < \SE(\mathcal{G})$ and $\LE(\mathcal{G}) = \SE(\mathcal{G})$. 
\end{prob} 
Das et al. \cite{DaAo-2018} also posed the following problem where  $\DE(\mathcal{G})$, $\DLE(\mathcal{G})$ and $\DQE(\mathcal{G})$ denote the distance energy, distance Laplacian energy and distance signless Laplacian energy of $\mathcal{G}$.
\begin{prob}\label{Prob 4}\cite[Problem 4]{DaAo-2018}
	Is there any connected graph $\mathcal{G}(\ncong K_n)$ such that $\E(\mathcal{G}) = \LE(\mathcal{G}) = \SE(\mathcal{G}) = \DE(\mathcal{G}) = \DLE(\mathcal{G}) = \DQE(\mathcal{G})$ ?
\end{prob}
Recently, Jannat and Nath \cite{JaNa7-2024} have computed several distance energies of NCCC-graphs of the groups considered in this paper and discussed relations (comparison) among them.  In Theorem \ref{Energy-NCCC1} and Corollary \ref{cor U6m}, we have seen that $\E(\mathcal{G}) = \LE(\mathcal{G}) = \SE(\mathcal{G})$.  In this section, we consider the groups $D_{2m}$, $T_{4m}$,  $SD_{8m}$, $V_{8m}$ and  compare  energy, Laplacian energy and signless Laplacian energy of NCCC-graphs of theses groups graphically.

\hspace{-.5cm}
\begin{minipage}[t]{.4\linewidth}
	\begin{tikzpicture}
		\begin{axis}
			[
			xlabel={$m$ is even and $\frac{m}{2}$ is odd ($m\geq 6$)  $\rightarrow$},
			ylabel={Energies of $\Gamma_{D_{2m}}$ $\rightarrow$},
			xmin=6, xmax=46,
			ymin=0, ymax=110,
			grid = both,
			minor tick num = 1,
			major grid style = {lightgray},
			minor grid style = {lightgray!25},
			width=.7\textwidth,
			height=.7\textwidth,
			legend style={legend pos=north west},
			% legend image post style={black}
			]
			\addplot[domain=6:46,samples at={6,10,14,18,22,26,30,34,38,42,46},mark=*,green, samples=24, mark size=.8pt]{2*sqrt(x-2)};
			\tiny
			\addlegendentry{$\E$}
			\addplot[domain=6:46,samples at={6,10,14,18,22,26,30,34,38,42,46},mark=triangle*,blue,mark size=.8pt, samples=24]{(2*(x-2)*(x-6))/(x+2) + 4};
			\tiny
			\addlegendentry{$\LE$}
			\addplot[domain=6:46,samples at={6,10,14,18,22,26,30,34,38,42,46},mark=square*, red, mark size=.5pt, samples=24]{(2*(x-2)*(x-6))/(x+2) + 4};
			\tiny
			\addlegendentry{$\SE$}
		\end{axis}
	\end{tikzpicture}
	\vspace{-.2 cm}
	\captionsetup{font=footnotesize}
	\captionof{figure}{Energies of $\Gamma_{D_{2m}}$, $m$ is even and $\frac{m}{2}$ is odd}\label{FigD2n1}
\end{minipage}
\hspace{-1cm}
\begin{minipage}[t]{.4\linewidth}
	\begin{tikzpicture}
		\begin{axis}
			[
			xlabel={$m$ and $\frac{m}{2}$ are even ($m\geq 4$) $\rightarrow$},
			ylabel={Energies of $\Gamma_{D_{2m}} \rightarrow$},
			xmin=4, xmax=36,
			ymin=0, ymax=70,
			grid = both,
			minor tick num = 1,
			major grid style = {lightgray},
			minor grid style = {lightgray!25},
			width=.7\textwidth,
			height=.7\textwidth,
			legend style={legend pos=north west},
			% legend image post style={black}
			]
			\addplot[domain=4:48,samples at={4,8,12,16,20,24,28,32,36,40,44,48},mark=*,green, samples=20, mark size=.8pt]{1 + sqrt(1 + 4*(x-2))};
			\tiny
			\addlegendentry{$\E$}
			\addplot[domain=4:48,samples at={4,8,12,16,20,24,28,32,36,40,44,48},mark=triangle*,blue,mark size=.8pt, samples=20]{(2*(x-2)*(x-4))/(x+2) + 4};
			\tiny
			\addlegendentry{$\LE$}
			\addplot[domain=8:48,samples at={8,12,16,20,24,28,32,36,40,44,48},mark=square*, red, mark size=.8pt, samples=20]{(2*(x-2)*(x-4))/(x+2) + ((x-2)*(sqrt(1+(16/(x-2)))-1))/2};
			\tiny
			\addlegendentry{$\SE$}
			\addplot[domain=4:4,samples at={4},mark=square*, red, mark size=.8pt, samples=20,forget plot]{4};
			\tiny
		\end{axis}
	\end{tikzpicture}
	\vspace{-.2 cm}
	\captionsetup{font=footnotesize}
	\captionof{figure}{Energies of $\Gamma_{D_{2m}}$, $m$ and $\frac{m}{2}$ are even}\label{FigD2n2}
\end{minipage}
\hspace{-1cm}
\begin{minipage}[t]{.4\linewidth}
	\begin{tikzpicture}
		\begin{axis}
			[
			xlabel={$m$ is odd ($m \geq 3$)  $\rightarrow$},
			ylabel={Energies of $\Gamma_{D_{2m}}$ $\rightarrow$},
			xmin=3, xmax=29,
			ymin=0, ymax=30,
			grid = both,
			minor tick num = 1,
			major grid style = {lightgray},
			minor grid style = {lightgray!25},
			width=.7\textwidth,
			height=.7\textwidth,
			legend style={legend pos=north west},
			% legend image post style={black}
			]
			\addplot[domain=3:29,samples at={3,5,7,9,11,13,15,17,19,21,23,25,27,29},mark=*,green, samples=24, mark size=.8pt]{2*(sqrt((x-1)/2))};
			\tiny
			\addlegendentry{$\E$}
			\addplot[domain=3:29,samples at={3,5,7,9,11,13,15,17,19,21,23,25,27,29},mark=triangle*,blue,mark size=.8pt, samples=24]{(((x-1)*(x-3))/(x+1)) + 2};
			\tiny
			\addlegendentry{$\LE$}
			\addplot[domain=3:29,samples at={3,5,7,9,11,13,15,17,19,21,23,25,27,29},mark=square*, red, mark size=.5pt, samples=24]{(((x-1)*(x-3))/(x+1)) + 2};
			\tiny
			\addlegendentry{$\SE$}
		\end{axis}
	\end{tikzpicture}
	\vspace{-.2 cm}
	\captionsetup{font=footnotesize}
	\captionof{figure}{Energies of $\Gamma_{D_{2m}}$, $m$ is odd }\label{FigD2n3}
\end{minipage}
%\hspace{0.05cm}

\begin{minipage}[t]{.5\linewidth}
	\begin{tikzpicture}
		\begin{axis}
			[
			xlabel={$m$ is even ($m \geq 2$)  $\rightarrow$},
			ylabel={Energies of $\Gamma_{T_{4m}}$ $\rightarrow$},
			xmin=2, xmax=30,
			ymin=0, ymax=110,
			grid = both,
			minor tick num = 1,
			major grid style = {lightgray},
			minor grid style = {lightgray!25},
			width=.7\textwidth,
			height=.7\textwidth,
			legend style={legend pos=north west},
			% legend image post style={black}
			]
			\addplot[domain=2:30,samples at={2,4,6,8,10,12,14,16,18,20,22,24,26,28,30},mark=*,green, samples=24, mark size=.8pt]{1 + sqrt(1 + 8*(x-1))};
			\tiny
			\addlegendentry{$\E$}
			\addplot[domain=2:30,samples at={2,4,6,8,10,12,14,16,18,20,22,24,26,28,30},mark=triangle*,blue,mark size=.8pt, samples=24]{((4*(x-1)*(x-2))/(x+1)) + 4};
			\tiny
			\addlegendentry{$\LE$}
			\addplot[domain=4:30,samples at={4,6,8,10,12,14,16,18,20,22,24,26,28,30},mark=square*, red, mark size=.5pt, samples=24]{((4*(x-1)*(x-2))/(x+1)) + (x - 1)*(sqrt(1 + (8/(x-1))) - 1)};
			\tiny
			\addlegendentry{$\SE$}
			\addplot[domain=2:2,samples at={2},mark=square*, red, mark size=.5pt, samples=20,forget plot]{4};
			\tiny
		\end{axis}
	\end{tikzpicture}
	\vspace{-.2 cm}
	\captionsetup{font=footnotesize}
	
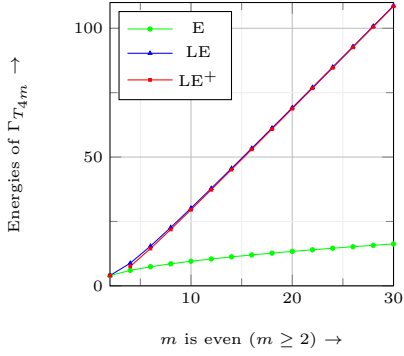
\captionof{figure}{Energies of $\Gamma_{T_{4m}}$, $m$ is even }\label{FigT4m1}
\end{minipage}
\hspace{0.05cm}
\begin{minipage}[t]{.5\linewidth}
	\begin{tikzpicture}
		\begin{axis}
			[
			xlabel={$m$ is odd ($m \geq 3$)  $\rightarrow$},
			ylabel={Energies of $\Gamma_{T_{4m}}$ $\rightarrow$},
			xmin=3, xmax=31,
			ymin=0, ymax=110,
			grid = both,
			minor tick num = 1,
			major grid style = {lightgray},
			minor grid style = {lightgray!25},
			width=.7\textwidth,
			height=.7\textwidth,
			legend style={legend pos=north west},
			% legend image post style={black}
			]
			\addplot[domain=3:31,samples at={3,5,7,9,11,13,15,17,19,21,23,25,27,29,31},mark=*,green, samples=24, mark size=.8pt]{4*(sqrt((x-1)/2))};
			\tiny
			\addlegendentry{$\E$}
			\addplot[domain=3:31,samples at={3,5,7,9,11,13,15,17,19,21,23,25,27,29,31},mark=triangle*,blue,mark size=.8pt, samples=24]{((4*(x-1)*(x-3))/(x+1)) + 4};
			\tiny
			\addlegendentry{$\LE$}
			\addplot[domain=3:31,samples at={3,5,7,9,11,13,15,17,19,21,23,25,27,29,31},mark=square*, red, mark size=.5pt, samples=24]{((4*(x-1)*(x-3))/(x+1)) + 4};
			\tiny
			\addlegendentry{$\SE$}
		\end{axis}
	\end{tikzpicture}
	\vspace{-.2 cm}
	\captionsetup{font=footnotesize}
	
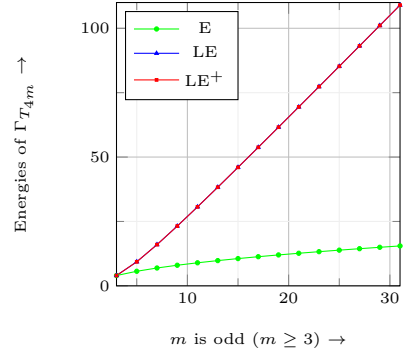
\captionof{figure}{Energies of $\Gamma_{T_{4m}}$, $m$ is odd }\label{FigT4m2}
\end{minipage}
\hspace{0.05cm}
%\begin{minipage}[t]{.5\linewidth}
%	\begin{tikzpicture}
%		\begin{axis}
%			[
%			xlabel={($m \geq 2$)  $\rightarrow$},
%			ylabel={Energies of $\Gamma_{U_{6m}}$ $\rightarrow$},
%			xmin=2, xmax=16,
%			ymin=0, ymax=40,
%			grid = both,
%			minor tick num = 1,
%			major grid style = {lightgray},
%			minor grid style = {lightgray!25},
%			width=.7\textwidth,
%			height=.7\textwidth,
%			legend style={legend pos=north west},
%			% legend image post style={black}
%			]
%			\addplot[domain=2:16,samples at={2,3,4,5,6,7,8,9,10,11,12,13,14,15,16},mark=*,green, samples=24, mark size=.8pt]{2*x};
%			\tiny
%			\addlegendentry{$\E$}
%			\addplot[domain=2:16,samples at={2,3,4,5,6,7,8,9,10,11,12,13,14,15,16},mark=triangle*,blue,mark size=.5pt, samples=24]{2*x};
%			\tiny
%			\addlegendentry{$\LE$}
%			\addplot[domain=2:16,samples at={2,3,4,5,6,7,8,9,10,11,12,13,14,15,16},mark=square*, red, mark size=.2pt, samples=24]{2*x};
%			\tiny
%			\addlegendentry{$\SE$}
%		\end{axis}
%	\end{tikzpicture}
%	\vspace{-.2 cm}
%	\captionsetup{font=footnotesize}
%	\captionof{figure}{Energies of $\Gamma_{U_{6m}}$ }\label{FigU6m}
%\end{minipage}
%\hspace{0.05cm}

\begin{minipage}[t]{.5\linewidth}
	\begin{tikzpicture}
		\begin{axis}
			[
			xlabel={$m$ is even ($m \geq 2$)  $\rightarrow$},
			ylabel={Energies of $\Gamma_{SD_{8m}}$ $\rightarrow$},
			xmin=2, xmax=20,
			ymin=0, ymax=160,
			grid = both,
			minor tick num = 1,
			major grid style = {lightgray},
			minor grid style = {lightgray!25},
			width=.7\textwidth,
			height=.7\textwidth,
			legend style={legend pos=north west},
			% legend image post style={black}
			]
			\addplot[domain=2:30,samples at={2,4,6,8,10,12,14,16,18,20,22,24,26,28,30},mark=*,green, samples=24, mark size=.8pt]{1 + (sqrt(1 + 8*(2*x-1)))};
			\tiny
			\addlegendentry{$\E$}
			\addplot[domain=2:30,samples at={2,4,6,8,10,12,14,16,18,20,22,24,26,28,30},mark=triangle*,blue,mark size=.8pt, samples=24]{((4*(2*x-1)*(2*x-2))/(2*x+1)) + 4};
			\tiny
			\addlegendentry{$\LE$}
			\addplot[domain=2:30,samples at={2,4,6,8,10,12,14,16,18,20,22,24,26,28,30},mark=square*, red, mark size=.4pt, samples=24]{((4*(2*x-1)*(2*x-2))/(2*x+1)) + (2*x - 1)*(sqrt(1 + (8/(2*x-1))) - 1)};
			\tiny
			\addlegendentry{$\SE$}
		\end{axis}
	\end{tikzpicture}
	\vspace{-.2 cm}
	\captionsetup{font=footnotesize}
	
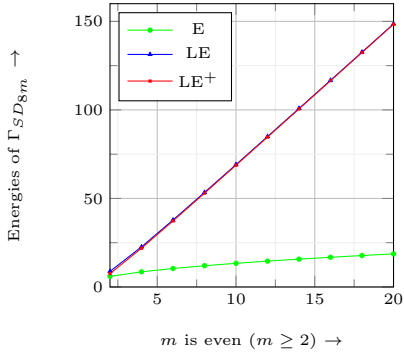
\captionof{figure}{Energies of $\Gamma_{SD_{8m}}$, $m$ is even }\label{FigSD8m1}
\end{minipage}
\hspace{0.05cm}
\begin{minipage}[t]{.5\linewidth}
	\begin{tikzpicture}
		\begin{axis}
			[
			xlabel={$m$ is odd ($m \geq 3$)  $\rightarrow$},
			ylabel={Energies of $\Gamma_{SD_{8m}}$ $\rightarrow$},
			xmin=3, xmax=21,
			ymin=0, ymax=300,
			grid = both,
			minor tick num = 1,
			major grid style = {lightgray},
			minor grid style = {lightgray!25},
			width=.7\textwidth,
			height=.7\textwidth,
			legend style={legend pos=north west},
			% legend image post style={black}
			]
			\addplot[domain=3:31,samples at={3,5,7,9,11,13,15,17,19,21,23,25,27,29,31},mark=*,green, samples=24, mark size=.8pt]{8*(sqrt((x-1)/2))};
			\tiny
			\addlegendentry{$\E$}
			\addplot[domain=3:31,samples at={3,5,7,9,11,13,15,17,19,21,23,25,27,29,31},mark=triangle*,blue,mark size=.8pt, samples=24]{((16*(x-1)*(x-3))/(x+1)) + 8};
			\tiny
			\addlegendentry{$\LE$}
			\addplot[domain=3:31,samples at={3,5,7,9,11,13,15,17,19,21,23,25,27,29,31},mark=square*, red, mark size=.5pt, samples=24]{((16*(x-1)*(x-3))/(x+1)) + 8};
			\tiny
			\addlegendentry{$\SE$}
		\end{axis}
	\end{tikzpicture}
	\vspace{-.2 cm}
	\captionsetup{font=footnotesize}
	
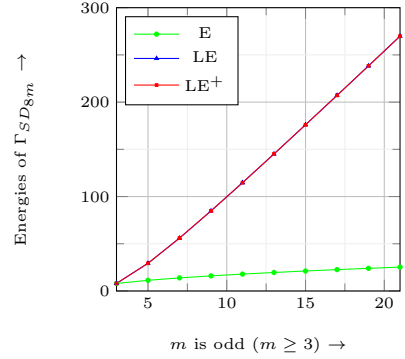
\captionof{figure}{Energies of $\Gamma_{SD_{8m}}$, $m$ is odd }\label{FigSD8m2}
\end{minipage}
%\hspace{0.05cm}

\begin{minipage}[t]{.5\linewidth}
	\begin{tikzpicture}
		\begin{axis}
			[
			xlabel={$m$ is even ($m \geq 2$)  $\rightarrow$},
			ylabel={Energies of $\Gamma_{V_{8m}}$ $\rightarrow$},
			xmin=2, xmax=20,
			ymin=0, ymax=300,
			grid = both,
			minor tick num = 1,
			major grid style = {lightgray},
			minor grid style = {lightgray!25},
			width=.7\textwidth,
			height=.7\textwidth,
			legend style={legend pos=north west},
			% legend image post style={black}
			]
			\addplot[domain=2:30,samples at={2,4,6,8,10,12,14,16,18,20,22,24,26,28,30},mark=*,green, samples=24, mark size=.8pt]{2*(1 + sqrt(1 + 8*(x-1)))};
			\tiny
			\addlegendentry{$\E$}
			\addplot[domain=2:30,samples at={2,4,6,8,10,12,14,16,18,20,22,24,26,28,30},mark=triangle*,blue,mark size=.8pt, samples=24]{((16*(x-1)*(x-2))/(x+1)) + 8};
			\tiny
			\addlegendentry{$\LE$}
			\addplot[domain=4:30,samples at={4,6,8,10,12,14,16,18,20,22,24,26,28,30},mark=square*, red, mark size=.5pt, samples=24]{(16*(x-1)*(x-2))/(x+1) + 2*(x - 1)*(sqrt(1 + 8/(x-1)) - 1)};
			\tiny
			\addlegendentry{$\SE$}
			\addplot[domain=2:2,samples at={2},mark=square*, red, mark size=.8pt, samples=20,forget plot]{8};
			\tiny
		\end{axis}
	\end{tikzpicture}
	\vspace{-.2 cm}
	\captionsetup{font=footnotesize}
	
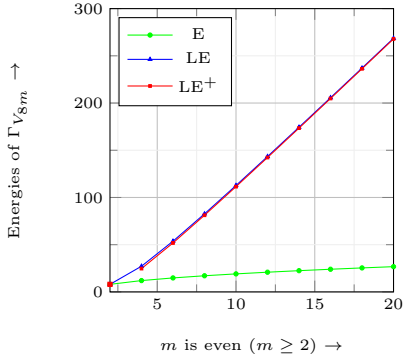
\captionof{figure}{Energies of $\Gamma_{V_{8m}}$, $m$ is even }\label{FigV8m1}
\end{minipage}
\hspace{0.05cm}
\begin{minipage}[t]{.5\linewidth}
	\begin{tikzpicture}
		\begin{axis}
			[
			xlabel={$m$ is odd ($m \geq 3$)  $\rightarrow$},
			ylabel={Energies of $\Gamma_{V_{8m}}$ $\rightarrow$},
			xmin=3, xmax=21,
			ymin=0, ymax=200,
			grid = both,
			minor tick num = 1,
			major grid style = {lightgray},
			minor grid style = {lightgray!25},
			width=.7\textwidth,
			height=.7\textwidth,
			legend style={legend pos=north west},
			% legend image post style={black}
			]
			\addplot[domain=3:31,samples at={3,5,7,9,11,13,15,17,19,21,23,25,27,29,31},mark=*,green, samples=24, mark size=.8pt]{1 + (sqrt(1 + 8*(2*x-1)))};
			\tiny
			\addlegendentry{$\E$}
			\addplot[domain=3:31,samples at={3,5,7,9,11,13,15,17,19,21,23,25,27,29,31},mark=triangle*,blue,mark size=.8pt, samples=24]{((4*x - 2)*(4*x - 2) + 8)/(2*x+1)};
			\tiny
			\addlegendentry{$\LE$}
			\addplot[domain=3:31,samples at={3,5,7,9,11,13,15,17,19,21,23,25,27,29,31},mark=square*, red, mark size=.5pt, samples=24]{((3*(2*x-1)*(2*x-3))/(2*x+1)) + (2*x - 1)*(sqrt(1 + (8/(2*x-1))))};
			\tiny
			\addlegendentry{$\SE$}
		\end{axis}
	\end{tikzpicture}
	\vspace{-.2 cm}
	\captionsetup{font=footnotesize}
	
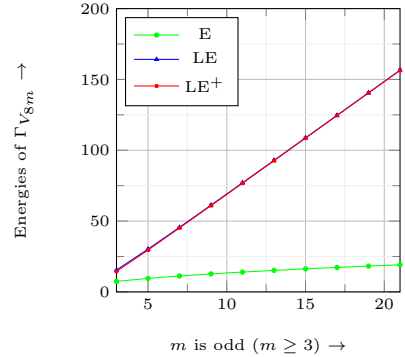
\captionof{figure}{Energies of $\Gamma_{V_{8m}}$, $m$ is odd }\label{FigV8m2}
\end{minipage}

We conclude this paper with the following result which is obtained from Figures \ref{FigD2n1} to \ref{FigV8m2}.
\begin{thm}
Let $G = D_{2m}, T_{4m}, SD_{8m}$ or $V_{8m}$. Then
\begin{enumerate}
	\item $\E(\Gamma_G) = \LE(\Gamma_G) = \SE(\Gamma_G)$ if and only if $G = D_6$, $D_8$, $D_{12}$, $T_8$, $T_{12}$, $SD_{24}$ or $V_{16}$.
	\item $\E(\Gamma_G) < \LE(\Gamma_G) = \SE(\Gamma_G)$ if and only if $G = D_{2m}$, where $m$ is odd or  $m \geq 10$ is even and $\frac{m}{2}$ is odd;  $T_{4m}$, where  $m \geq 5$ is odd; $SD_{8m}$, where $m \geq 5$ is odd.
	\item $\E(\Gamma_G) < \SE(\Gamma_G) < \LE(\Gamma_G)$ if and only if $G = D_{2m}$, where $m \geq 8$ and $\frac{m}{2}$ is even; $T_{4m}$, where  $m \geq 4$ is even; $SD_{8m}$, where  $m$ is even; $V_{8m}$, where $m \geq 3$.    
	\end{enumerate}
\end{thm}

\end{document}